\documentclass[11pt]{article}

\usepackage{amsfonts}
\usepackage{amssymb}
\usepackage{amsmath}

\usepackage{amsthm}
\usepackage{algorithm}
\usepackage{enumerate}
\usepackage{tabu}
\usepackage{natbib}  

\usepackage[colorlinks,linkcolor=blue,anchorcolor=green,citecolor=black]{hyperref}

=0pt
\def\sqr#1#2{{\vcenter{\vbox{\hrule height.#2pt
              \hbox{\vrule width.#2pt height#1pt \kern#1pt \vrule
width.#2pt}
              \hrule height.#2pt}}}}
\def\dbE{{\mathbb{E}}}
\def\dbF{{\mathbb{F}}}

\def\dbH{{\mathbb{H}}}

\def\dbN{{\mathbb{N}}}

\def\dbP{{\mathbb{P}}}

\def\dbR{{\mathbb{R}}}
\def\dbS{{\mathbb{S}}}

\def\dbU{{\mathbb{U}}}
\def\dbV{{\mathbb{V}}}
\def\dbW{{\mathbb{W}}}
\def\dbX{{\mathbb{X}}}

\def\a{\alpha}
\def\b{\beta}
\def\g{\gamma}

\def\l{\lambda}

\def\si{\sigma}
\def\t{\tau}
\def\f{\varphi}

\def\om{\omega}

\def\3n{\negthinspace \negthinspace \negthinspace }
\def\2n{\negthinspace \negthinspace }
\def\1n{\negthinspace }
\def\ns{\noalign{\smallskip} }

\def\ds{\displaystyle}
\def\G{\Gamma}
\def\D{\Delta}

\def\L{\Lambda}

\def\Om{\Omega}
\def\om{\omega}
\def\cC{{\cal C}}

\def\cF{{\cal F}}
\def\cG{{\cal G}}
\def\cH{{\cal H}}
\def\cI{{\cal I}}
\def\cJ{{\cal J}}

\def\cL{{\cal L}}

\def\cX{{\cal X}}

\def\mE{{\mathbb{E}}}

\def\no{\noindent}

\def\ss{\smallskip}
\def\ms{\medskip}

\def\q{\quad}
\def\qq{\qquad}

\def\pa{\partial}
\def\h{\widehat}
\def\wt{\widetilde}
\def\cd{\cdot}
\def\cds{\cdots}

\def\ae{\hbox{\rm a.e.{ }}}
\def\as{\hbox{\rm a.s.{ }}}

\def\span{\hbox{\rm span$\,$}}

\def\deq{\mathop{\buildrel\D\over=}}

\def\({\Big (}
\def\){\Big )}
\def\[{\Big[}
\def\]{\Big]}

\def\={\buildrel \triangle \over =}

\def\resp{{\it resp. }}

\def\lt{\left}
\def\rt{\right}
\def\bal{\begin{aligned}}
\def\eal{\end{aligned}}
\def\beq{\begin{equation*}\begin{aligned}}
\def\eeq{\end{aligned}\end{equation*}}
\def\O{\Omega}
\def\mf{\mathcal{F}}
\def\me{\mathbb{E}}
\def\llan{\left\langle}
\def\rran{\right\rangle}
\def\lan{\langle}
\def\ran{\rangle}
\def\rf{\eqref}

\def\be{\begin{equation}}
\def\bel{\begin{equation}\label}
\def\ee{\end{equation}}
\def\bea{\begin{eqnarray}}
\def\eea{\end{eqnarray}}
\def\bt{\begin{theorem}\label}
\def\et{\end{theorem}}
\def\bc{\begin{corollary}\label}
\def\ec{\end{corollary}}
\def\bl{\begin{lemma}\label}
\def\el{\end{lemma}}
\def\bp{\begin{proposition}\label}
\def\ep{\end{proposition}}
\def\br{\begin{remark}\label}
\def\er{\end{remark}}
\def\ba{\begin{array}}
\def\ea{\end{array}}
\def\bd{\begin{definition}\label}
\def\ed{\end{definition}}

\newtheorem{lemma}{Lemma}[section]
\newtheorem{remark}{Remark}[section]

\newtheorem{theorem}{Theorem}[section]
\newtheorem{corollary}{Corollary}[section]

\newtheorem{definition}{Definition}[section]
\newtheorem{proposition}{Proposition}[section]

\allowdisplaybreaks
 \makeatletter
   
   \@addtoreset{equation}{section}
\makeatother

\begin{document}

\title{\bf Numerics for Stochastic Distributed Parameter Control Systems: a Finite Transposition
Method}

\author{Qi L\"u,\thanks{
School of Mathematics, Sichuan
University, Chengdu 610064, China. The
research of this author is supported by
the NSF of China under grants 12025105, 11971334
and 11931011, and  the Chang Jiang
Scholars Program from the Chinese
Education Ministry. {\small\it E-mail:}
{\small\tt lu@scu.edu.cn}. \ms}~~~
Penghui Wang,\thanks{ School of
Mathematics, Shandong University, Jinan
250100, China. The research of this
author is supported by the NSF of China
under grant 11871308. {\small\it
E-mail:} {\small\tt phwang@sdu.edu.cn}.
\ms} ~~~Yanqing Wang\thanks{ School of
Mathematics and Statistics, Southwest
University, Chongqing 400715, China.
The research of this author is
supported by the NSF of China under
grant 11801467, and the Chongqing
Natural Science Foundation under grant
cstc2018jcyjAX0148. {\small\it E-mail:}
{\small\tt yqwang@amss.ac.cn}. \ms}~~~
  and~~~
Xu Zhang\thanks{School of Mathematics,
Sichuan University, Chengdu 610064,
 China. The research
of this author is partially supported by the NSF
of China under grants 11931011 and 11821001. {\small\it E-mail:}
{\small\tt zhang$\_$xu@scu.edu.cn}.}}

\date{}

\maketitle

\begin{abstract}
In this chapter, we present some recent progresses on the
numerics for stochastic distributed parameter control systems, based on the \emph{finite transposition
method} introduced in our previous works. We first explain how to reduce the
numerics of some stochastic control problems in this respect to the numerics of backward
stochastic evolution equations. Then we present a method to find finite
transposition solutions to such equations. At last, we give an illuminating example.  
\end{abstract}

\vspace{+0.3em}


\vspace{+0.3em}



\section{Introduction}


In this chapter, we study the numerics for stochastic distributed parameter control systems, including particularly
controlled stochastic evolution equations in infinite dimensions. Although some
pioneer works appeared in 1960s \citep[e.g.,][]{Bensoussan69, Kushner68, Tzafestas-Nightingale68},
 the  theory for
stochastic distributed parameter control systems is far from
mature till now. Lots of interesting and
important problems therein are not well studied. One of
them is the numerical computation on 
control problems of such systems, for which one will meet many
substantially difficulties:
\begin{itemize}

\item For stochastic control systems,
   numerical schemes should keep the adaptedness
   of states and controls with respect to the
   filtration. This restricts the application of
   many implicit schemes.

\item Unlike the deterministic setting, solutions to 
stochastic evolution equations (even in finite
dimensions) are usually non-differentiable with respect to
the time variable. This leads to some serious troubles
in the analysis of the convergence rate.

\item In the continuous-time setting, similar to the deterministic situation, a typical way to solve many stochastic control problems is to introduce the corresponding dual systems, which are however backward
stochastic evolution equations usually \citep[see][for more details]{Lv-Zhang20}. 
Numerically, it is
very hard to compute the correction terms in such sort of backward equations
since there is
no efficient method to compute the conditional
expectation of a general random variable.
\end{itemize}
Due to these difficulties, there are very few
works on the  numerical analysis for 
 stochastic distributed
parameter control systems
\citep[e.g.,][]{Dunst-Prohl16, Prohl-Wang20-2, Prohl-Wang20, Li-Zhou20}.
 
Because of the importance of backward
stochastic evolution equations in the study of stochastic distributed
parameter control systems, in this chapter we shall focus on the numerical solutions to such equations and the related theoretical analysis.

Up to now, there exist quite a number
of numerical schemes for backward
stochastic differential
equations, i.e. backward
stochastic evolution equations in finite dimensions, such as 
the method of four-step scheme \citep[e.g.,][]{Douglas-Ma-Protter96}, the
quantization tree method \citep[e.g.,][]{Bally-Pages03}, 
the forward/backward Euler method \citep[e.g.,][]{ZhangJF04,Gobet-Lemor-Warin05,Bender-Denk07}, the Malliavin
calculus based method \citep[e.g.,][]{Bouchard-Touzi04, Hu-Nualart-Song11}, the fully time-space discretization method
\citep[e.g.,][]{Zhao-Chen-Peng06, Zhao-Fu-Zhou14}, the finite
transposition method \citep[e.g.,][]{Wang-Zhang11}, the Wiener
chaos decomposition method \citep[e.g.,][]{Briand-Labart14}, the machine
learning method \citep[e.g.,][]{E-Hutzenthaler-Jentzen-Kruse19},  etc. However,
numerical schemes for backward
stochastic evolution equations in infinite dimensions (including typically backward
stochastic partial differential
equations) are quite
limited \citep{WangY13, Wang16, Dunst-Prohl16}.

Based on the stochastic transposition
method introduced in
\cite{Lv-Zhang13, Lv-Zhang14}, in this chapter, we propose a new numerical
algorithm to solve backward
stochastic evolution equations, which can be
used to compute the desired control of some
controllability/optimal control problems for stochastic distributed
parameter systems numerically. 

To present the key idea in the simplest
way, we do not pursue the full
technical generality. Firstly, we
consider only the simplest case of one
dimensional standard Brownian motion
(with respect to the time variable $t$). It
would not be difficult to extend the
results to the case of the $Q$-Brownian motion  and the
cylindrical Brownian motion. Secondly,
we impose considerably strong
regularity and boundedness assumptions
on data 
appeared in the
state equations and the
cost functionals.  Thirdly, we
assume that the control operators are bounded linear
ones, which cannot cover the case
of stochastic partial differential
equations with boundary/pointwise
controls. All of these
assumptions can be considerably relaxed. Fourthly, we only
consider a null controllability and a
linear quadratic optimal control problem. More
general controllability and optimal
control problems can be studied by the
method presented in this chapter.

The rest of this chapter is organized as follows. In Section \ref{sec2}, we shall introduce the control problems and the corresponding backward stochastic evolution equations to be considered in this chapter.  
In Section \ref{FTS-1}, we firstly
give an outline of the finite
transposition method for solving backward stochastic evolution equations numerically, and then provide
a choice of finite transposition space.
In Section \ref{FTS-2}, based on the
finite transposition space, we define
the finite transposition solution to
backward stochastic evolution equations, and
then present the existence and uniqueness of
this solution as well as the convergence rate of this numerical method.
Finally, in Section \ref{Ap}, as an application, we propose a numerical scheme
for an optimal control problem, based on our finite transposition
method.

\section{Adjoint equations for stochastic distributed parameter control problems}\label{sec2}

We begin with some notations and notions
to be used later.
Let  $T>0$ and $(\Omega,\cF,\mathbf{F},\dbP)$
(with $\mathbf{F}=\{\cF_t\}_{t\in [0,T]}$) be a
 filtered probability space so that $(\Omega,\cF,\dbP)$ is
complete, $\cF_0$ contains all $\dbP$-null sets
in $\cF$, and $\mathbf{F}$ is right continuous.
Let $\{W(t)\}_{t\in[0,T]}$ be an
$\mathbf{F}$-adapted $1$-dimensional standard
Brownian motion on $(\Om,\cF,\mathbf{F},\dbP)$.
Denote by $\dbF$ the progressive
$\si$-field (in $[0,T]\times\Omega$)
with respect to $\mathbf{F}$, by $\dbE
\xi$ the (mathematical) expectation of an
integrable random variable
$\xi:(\Omega,\cF,\dbP)\to \dbR$, and by
$\cC$ a generic positive constant,
which may vary from one place
to another.

For a Banach space $\cX$, we write $\|\cd\|_\cX$ for its norm. If $\cX$ is a Hilbert space, we denote its inner product by $\lan \cd,\cd\ran_\cX$ (When $\cX$ is a
Euclidean space, for simplicity, we denote its norm and inner product by $|\cd|$ and $\lan \cd,\cd\ran$, respectively).
For any filtration $\mathbf{G}=\{\cG_t\}_{t\in[0,T]}\subseteq \mathbf{F}$,
$p,q\in[1,\infty)$ and $t\in [0,T]$, let
$L_{\cG_t}^p(\Omega;\cX)\deq
L^p(\Omega,\cG_t,\dbP;\cX)$, and
$$
\begin{array}{ll}
\ns\ds L^q_\dbF(0,T;L^p(\Omega;\cX)) \deq \Big\{
\f :(0,T) \times \Omega\!\to\cX \,\bigm|\,
\f(\cd)\hbox{
is $\mathbf{F}$-adapted}\\
\ns\ds\qq\qq\qq\qq\qq \mbox{  and
}\int_0^T\(\dbE\|\f(t)\|_\cX^p\)^{\frac{q}
{p}}dt<\infty\Big\},\\
\ns\ds L^{p}_{\dbF}(\Omega;C([0,T];\cX))\deq
\Big\{\f:[0,T]\times\Omega\to \cX \,\bigm|\,
\f(\cd)\hbox{
is continuous, }\\
\ns\ds\qq\qq\qq\qq\qq \mathbf{F}\mbox{-adapted and
}\dbE\big(\|\f(\cdot)\|_{C([0,T];\cX)}^p\big)<\infty\Big\},\\
\ns\ds C_{\dbF}([0,T];L^{p}(\Omega;\cX))\deq
\Big\{\f:[0,T]\times\Omega\to
\cX\,\bigm|\,\f(\cd)\hbox{
is $\mathbf{F}$-adapted }\\
\ns\ds\qq\qq\qq\qq\qq \mbox{ and }\f(\cd):[0,T]
\to L^p_{\cF_T}(\Omega;\cX)\mbox{ is
continuous}\Big\}.
\end{array}
$$
In a similar way, we can define 
$
L^\infty_{\cG_t}(\O;\cX)$ and $ L^\infty_\dbF(0,T;L^\infty(\Omega;\cX)).
$
One can show that, for $1\le p,q\le \infty$,
$
L^p_{\cG_t}(\O;\cX)$ and $L^q_\dbF(0,T;L^p(\Omega;\cX))$ are Banach spaces (with the canonical norms). 
Further, both $L^{p}_{\dbF}(\Omega;C([0,T];\cX))$ and
$C_{\dbF}([0,T];L^{p}(\Omega;\cX))$ are 
 Banach spaces with norms given by  
$\|\f(\cd)\|_{L^{p}_{\dbF}(\Omega;C([0,T];\cX))}=\big(\dbE(\|\f(\cdot)\|_{C([0,T];\cX)}^p)\!\big)^{1/p}$ and
$\|\f(\cd)\|_{C_{\dbF}([0,T];L^{p}(\Omega;\cX))}=\max_{t\in
[0,T]}\big(\dbE(\|\f(t)\|_\cX^p)\big)^{1/p}$,
respectively.  
Also, we write
$D_{\dbF}([0,T];L^{p}(\Omega;\cX))$ for the
Banach space of all
$L^{p}(\Omega;\cX)$-valued,
$\mathbf{F}$-adapted, c\`adl\`ag 
stochastic processes $\f(\cd)$ such
that
$\max_{t\in
[0,T]}\big(\dbE(\|\f(t)\|_\cX^p)\big)^{1/p}
< \infty$, with the canonical norm. 
In the
sequel, we shall simply denote
$L^p_\dbF(0,T;L^p(\Omega;\cX))$ by
$L^p_\dbF(0,T;\cX)$; and further simply denote $L^p_{\cG_t}(\O;\dbR)$,
$L^p_\dbF(0,T;\dbR)$ and $L^p_{\dbF}(\Omega; C([0,T];\dbR))$ by $L^p_{\cG_t}(\O)$, $L^p_\dbF(0,T)$
and $L^p_{\dbF}(\Omega;C([0,T]))$, respectively.


In what follows, we assume that $H$
and $U$ are (real) separable Hilbert spaces, and that $A$ is
an unbounded linear operator (with domain $D(A)\subseteq H$), which generates a $C_0$-semigroup
$\{e^{At}\}_{t\geq 0}$ on $H$. Write $A^*$ for the
adjoint operator of $A$.
For any $y_0\in H$,  consider the following
controlled linear stochastic evolution equation:
\bel{system1}
\left\{
\bal
& dy(t) =\big(Ay(t) + B u(t)
\big)dt + \big(C y(t) +D u(t) \big)dW(t), \q t\in (0,T],\\
& y(0)=y_0.
\eal
\right.
\ee
%
In \eqref{system1}, $C(\cd)\in
L^\infty_\dbF(0,T;\cL(H))$, $B(\cd), D(\cd)\in
L^\infty_\dbF(0,T;\cL(U;H))$, $ u(\cd)\in
 L^2_\dbF(0,T;U)$ is the control
variable, and $y(\cd)=y(\cd;y_0,u(\cd))$ is the
state variable. By the well-posedness result for
stochastic evolution equations 
\citep[e.g.,][Section 3.2]{Lv-Zhang20},   the system
\eqref{system1} admits a unique mild solution
$y(\cd) \in C_\dbF([0,T];L^2(\Om;H))$.

\ms

Now, we recall the notion of null
controllability of \eqref{system1}.
\begin{definition}\label{def-exact-con}
System \eqref{system1} is called  null
controllable at time $T$ if for any $y_0\in H$,
there exists $u(\cd)\in L^2_\dbF(0,T;U)$
such that the
corresponding mild solution to \eqref{system1}
satisfies $y(T)=0$, a.s.
\end{definition}

Next, we introduce 
the following quadratic cost functional
\begin{equation}\label{cost}
\begin{array}{ll}\ds \cJ(y_0;u(\cd))
=\frac{1}{2}\mE\Big[ \int_0^T \big(\big\langle Q
y(t),y(t)\big\rangle_H +\big\langle R
u(t),u(t)\big\rangle_U\big)dt + \langle
Gy(T),y(T)\rangle_H\Big],
\end{array}
\end{equation}
where $Q(\cd)\in L^\infty_\dbF(0,T;\cL(H))$ and $R(\cd)\in
L^\infty_\dbF(0,T;\cL(U))$  are self-adjoint
operator-valued stochastic processes, while
$G\in L^\infty_{\cF_T}(\Om;\cL(H))$ is a
self-adjoint operator-valued random variable. 
Consider the linear quadratic
optimal control problem (SLQ problem
for short) as follows:

\ms

\no\bf Problem (SLQ): \rm For each $ y_0 \in
H$, find a $ \bar u(\cd)\in L^2_\dbF(0,T;U)$ such
that
\begin{equation}\label{5.2-eq3}
\cJ\big(y_0; \bar u(\cd)\big)=\inf_{u(\cd)\in
L^2_\dbF(0,T;U)}\cJ\big(y_0;u(\cd)\big).
\end{equation}
Any $ \bar u(\cdot)$ satisfying \rf{5.2-eq3}
is called an {\it optimal control}, the
corresponding state $ \bar y(\cdot)$ is called an
{\it optimal state}, and $\big(\bar y(\cdot), \bar u(\cdot)\big)$ is called an {\it
optimal pair}.

As in the deterministic case, in order to
solve the above null controllability
and optimal control problems, one may
employ the duality argument. For this
purpose, people introduce respectively the
following two backward stochastic
evolution equations:
\begin{equation}\label{bsystem1}
\left\{
\bal
&dz(t)
=-\big(A^*z(t) +C^*Z(t) \big)dt +Z(t)
dW(t), \q t\in [0,T),\\
& z(T)= z_T (\in L^2_{\cF_T}(\Om;H)),
\eal
\right.
\end{equation}
and
\begin{equation}\label{bsystem2}
\left\{
\bal
&dz(t)=-\big(A^*z(t)-Q y(t)+C^*Z(t)  \big)dt +Z(t)dW(t), \q t\in [0,T),\\
& z(T)= -Gy(T).
\eal
\right.
\end{equation}
These two equations can be used to solve aforementioned two stochastic control problems respectively (see Theorem \ref{th-nu-va} and Theorem \ref{LQth max}).

As far as we know, the study of backward stochastic
evolution equations is stimulated by
the works
\cite{Bensoussan83, Hu-Peng91}. Now, backward stochastic
evolution equations and its
variants play fundamental
roles in the theory of stochastic distributed parameter
control systems
\citep[e.g.,][]{Lv-Zhang20}.

Since neither the usual
natural filtration condition nor the
quasi-left continuity is assumed for
the filtration $\mathbf{F}$, one cannot
apply the existing results on infinite
dimensional backward stochastic
evolution equations 
\citep[e.g.,][]{Hu-Peng91, Al-Hussein09} to
obtain the well-posedness of the equations \eqref{bsystem1} and
\eqref{bsystem2} in
the sense of mild or weak solutions. Therefore, we shall   employ the
stochastic transposition method,
developed first in our paper
\cite{Lv-Zhang13} for backward stochastic differential equations, then
in \cite{Lv-Zhang14}  for backward
stochastic evolution equations, to
study the well-posedness of the equations \eqref{bsystem1} and
\eqref{bsystem2}.

Consider the following backward stochastic
evolution equation (in a general form):
\begin{equation}\label{bsystem3} \left\{ \bal &dz(t) =
-A^*z(t)dt + F(t,z(t),Z(t))dt + Z(t)
dW(t), \q t\in [0,T),\\
& z(T)=z_T.
\eal
\right.
\end{equation}
%
Here  $z_T\in L^2_{\cF_T}(\Om;H)$  and
$F:[0,T]\times\Omega\times H\times H\to H$ is a
given function  satisfying that
\bel{c3-Lip}
\left\{
\bal
&F(\cd,0,0)\in L^1_{\dbF}(0,T;L^2(\Omega;H)), \\
&  \|F(t,y_1,z_1)-F(t,y_2,z_2)\|_{H} \leq
\cC(\|y_1-y_2\|_H + \|z_1-z_2\|_{H}),
\\
& \qq\qq\qq\qq\forall\;
y_1,y_2,z_1,z_2\in {H}, \ \ \ae (t,\om)\in
[0,T]\times\Om.
\eal
\right.
\ee
%

To define the transposition solution to
\eqref{bsystem3}, we introduce the following
stochastic evolution equation:
\bel{fsystem2}
\left\{
\bal
& d\f(s) = (A\f(s) + v_1(s))ds +  v_2(s) dW(s), \q s\in (t,T],\\
& \f(t)=\eta,
\eal
\right.
\ee
%
where $t\in[0,T]$, $v_1(\cd)\in
L^1_{\dbF}(t,T;L^{2}(\Omega;H))$,
$v_2(\cd)\in L^{2}_{\dbF}(t,T;H)$ and
$\eta\in L^{2}_{\cF_t}(\Omega;H)$. By
the classical well-posedness result for
stochastic evolution equations 
\citep[e.g.,][Section 3.2]{Lv-Zhang20}, the equation \eqref{fsystem2}
admits a unique mild solution
$\f(\cd)\in
C_\dbF([t,T];L^2(\Omega;H))$, and
\beq
\bal
\|\f(\cd)\|_{C_\dbF([t,T];L^2(\Omega;H))}
\leq \cC\big(\|\eta\|_{L^{2}_{\cF_t}(\Omega;H)} +
\|v_1(\cd)\|_{L^1_{\dbF}(t,T;L^{2}(\Omega;H))} +
\|v_2(\cd)\|_{L^{2}_{\dbF}(t,T;H)} \big).
\eal
\eeq

If \eqref{bsystem3} admits a classical
mild solution in
\cite{Hu-Peng91, Al-Hussein09} (say,
when $\mathbf{F}$ is the natural filtration of
the Brownian motion $W(\cd)$), then by
It\^o's formula, we have
\begin{equation}\label{eq def solzz}
\begin{array}{ll}\ds
\dbE \big\langle
\f(T),z_T\big\rangle_{H}
-\dbE\int_t^T \big\langle \f(s),F(s,z(s),Z(s) )\big\rangle_Hds\\
\ns\ds = \dbE
\big\langle\eta,z(t)\big\rangle_H +
\dbE\int_t^T \big\langle
v_1(s),z(s)\big\rangle_H ds +
\dbE\int_t^T \big\langle
v_2(s),Z(s)\big\rangle_{H} ds.
\end{array}
\end{equation}
Motivated by this, we introduce the
following notion.
\begin{definition}\label{definition1}
We call $ (z(\cdot), Z(\cdot)) \in
D_{\dbF}([0,T];L^{2}(\Omega;H)) \times
L^2_{\dbF}(0,T;H)$ a transposition
solution to the equation
\rf{bsystem3} if for any $t\in
[0,T]$, $v_1(\cdot)\in
L^1_{\dbF}(t,T;L^2(\Omega;H))$,
$v_2(\cdot)\in L^2_{\dbF}(t,T;H)$,
$\eta\in L^2_{\cF_t}(\Omega;H)$ and the
corresponding solution $\f(\cd)\in
C_{\dbF}([t,T];L^2(\Omega;H))$ to
\eqref{fsystem2}, the equality
\eqref{eq def solzz} holds.
\end{definition}

The space for the
first component of the solution is
chosen to be $D_\dbF([0,
T];L^2(\Om;H))$ rather than $C_\dbF([0,
T];L^2(\Om;H))$. This is quite natural
because the filtration $\mathbf{F}$ is
assumed only to be right-continuous.

It is easy to see that if
\eqref{bsystem3} has a mild solution,
then such a solution is also a
transposition solution to \eqref{bsystem3}.

\begin{remark}
The above stochastic transposition method  is stimulated by the
classical transposition method to solve the
non-homogeneous boundary value problems for
deterministic partial differential equations \citep[see, e.g.,][]{Lions-Magenes72}.
This method is a variant of the
standard duality method, and in some
sense it provides a way to see
something which is not easy to be
detected directly. Specifically, for
the equation \eqref{bsystem3}, the
point of this method is to interpret
the solution to a backward stochastic
evolution equation in terms of a
forward stochastic evolution equation
which is well studied. The key tool to do this is the Riesz type representation theorem for $L^q_\dbF(0,T;L^p(\Omega;\cX))$ $(1\leq p,q<\infty)$ \citep[see][]{Lu-Yong-Zhang12, Lu-Yong-Zhang18}.  On the other
hand, the equality \eqref{eq def solzz}
can be regarded as a variational formulation
of the equation \eqref{bsystem3}. As
pointed out in \citep[Remark
3.2]{Lv-Zhang13}, it provides a way to
solve the equation \eqref{bsystem3}
numerically.
\end{remark}

We have the following well-posedness result for
the equation \rf{bsystem3}.
\begin{theorem}[{\citet[][Theorem 3.1]{Lv-Zhang14}}]\label{th1} 
For any $z_T \in L^2_{\cF_T}(\Omega;H)$, the
equation \eqref{bsystem3} admits a unique
transposition solution $(z(\cdot), Z(\cdot)) \in
D_{\dbF}([0,T];L^{2}(\Omega; H)) \times
L^2_{\dbF}(0,T;H)$. Furthermore,
\begin{equation*}\label{th1-eq1}
\bal
\|(z(\cdot), Z(\cdot))\|_{
D_{\dbF}([0,T];L^2(\Omega;H)) \times
L^2_{\dbF}(0,T;H)}\leq \cC\big(\|z_T\|_{
L^2_{\cF_T}(\Omega;H)}+ \|F(\cdot,0,0)\|_{
L^1_{\dbF}(0,T;L^2(\Omega;H))}\big).
\eal
\end{equation*}
\end{theorem}

From Theorem \ref{th1}, it follows that both
\eqref{bsystem1} and \eqref{bsystem2} are
well-posed in the sense of transposition
solution. Then we can study the null
controllability and optimal control problems by
\eqref{bsystem1} and \eqref{bsystem2},
respectively.

We first consider the problem of finding the
control which drives the state of the system
\eqref{system1} to rest. To this end, define
a functional on $L^2_{\cF_T}(\Om;H)$ as follows:
\bel{th-ex-va-eq2}
\bal
\mathfrak{J}(z_T)  = \frac{1}{2}
\mE\int^T_0 \|B^* z(t) +D^* Z(t)\|^2_{U} dt  +\mE\langle y_0,z(0)\rangle_{H},\q \forall\,z_T\in
L^2_{\cF_T}(\Om;H),
\eal
\ee
where $(z(\cd),Z(\cd))=(z(\cd;z_T),Z(\cd;z_T))$ is the transposition
solution to the equation
\eqref{bsystem1} (corresponding to the final datum $z_T$). 

The following result holds:
\begin{theorem}[{\citet[Theorems 7.17 and 7.28]{Lv-Zhang20}}]\label{th-nu-va} 
If the system \eqref{system1} is null
controllable at time $T$, then, among all
controls transferring the state of
\eqref{system1} from $y_0$ to $0$ at time $T$,
the one given by
\bel{th-ex-va-eq1}
\bal
u(t)=  B^*  z(t;\hat z_T)+ D^* Z(t;\hat z_T), \q\ae
(t,\om)\in[0,T]\times\Om 
\eal
\ee
has the minimal $L^2_\dbF(0,T;U)$-norm, where $\hat
z_T\in L^2_{\cF_T}(\Om;H)$ is the minimizer of
the functional $\mathfrak{J}(\cd)$.
\end{theorem}
\begin{remark}
	In this chapter, since we focus on the
	numerics, we will not consider problem that when \eqref{system1} is null controllable, which can be reduced to
	an observability estimate of the
	equation \eqref{bsystem1} 
	\citep[e.g.,][Theorem 7.17]{Lv-Zhang20}.
\end{remark}

Theorem \ref{th-nu-va} provides an explicit
formula of the desired control for null controllability of \eqref{system1}. Therefore,
to get such a control, one has to solve the
following two problems:
\begin{enumerate}[(1)]
 \item Finding the minimizer $\hat
z_T$ of the functional $\mathfrak{J}(\cd)$;

  \item Solving \eqref{bsystem1} to obtain $B^*  z(t;\hat z_T)+ D^* Z(t;\hat
  z_T)$.
\end{enumerate}
If we know how to solve the equation
\eqref{bsystem1}, then we can adopt the
gradient method/Newton's method  to
find the minimizer of the functional
$\mathfrak{J}(\cd)$ numerically. Hence, the key
point to compute the control \eqref{th-ex-va-eq1} is
how to solve the equation \eqref{bsystem1}.

\ms

Next, we consider  Problem (SLQ). With the aid
of \eqref{bsystem2}, we have the following
result.
\begin{theorem}[{\citet[Theorem 5.2]{Lv-Zhang15}}]\label{LQth max} 
Let $(\bar y(\cd), \bar u(\cd))$ be  an optimal
pair of Problem (SLQ). Then, for the
transposition solution $(z(\cd), Z(\cd))$ to
\eqref{bsystem2} with $y(\cd)$ replaced by $\bar y(\cd)$,
it holds that
\begin{equation}\label{max pr1}
R(t) \bar u(t)-B(t)^* z(t)-D(t)^*Z(t) =0, \q\ae
(t,\om)\in[0,T]\times\Om.
\end{equation}
\end{theorem}
According to Theorem \ref{LQth max}, by \eqref{system1} and \rf{bsystem2}, when $R(\cd)^{-1}\in L^\infty_\dbF(0,T;\cL(U))$, in order for finding the optimal
pair $(\bar y(\cd), \bar u(\cd))$ of Problem (SLQ), it suffices to solve the following coupled forward-backward stochastic
evolution equation:
\begin{equation}\label{LQfbsystem3}
\left\{\!\!\!
\begin{array}{ll}\ds
d\bar y(t)=\big(A \bar y(t)+ B R^{-1} B^* z(t) + B
R^{-1} D^*  Z(t) \big)dt\\
\ns\ds\hspace{1.1cm} +\big( C \bar y(t) + D R^{-1}
B^* z(t) + D R^{-1} D^* Z(t)\big)dW(t)  &\q  t\in (0,T],
\\
\ns\ds dz(t)=- \big(A^* z(t) -Q\bar y(t) + C^* Z(t)
\big)dt+  Z(t) dW(t) &\q t\in [0,T),\\
\ns\ds \bar y(0)=y_0,\qq z(T)=-G\bar y(T).
\end{array}
\right.
\end{equation}
Nevertheless, a key
point to solve \rf{LQfbsystem3} is to handle the backward equation (for $(z(\cd),Z(\cd))$) therein.

From the above discussion, we see that it is very
important to solve backward stochastic evolution
equations in the study of  control problems
for stochastic evolution equations.
Note that, for almost all backward stochastic
evolution equations, it is impossible to find an
explicit formula of the solutions. Hence, for practical applications, it is
crucial to find numerical solutions to these
equations.

\section{The space of finite transposition}\label{FTS-1}

In the following two sections, based on
the  transposition method, we shall present the
{\em
finite transposition method}, introduced in \citep{Wang-Zhang11, Wang-Wang-Zhang20} for backward stochastic differential equations,  to find
numerical solutions for backward
stochastic evolution equations. To
present this method clearly, we only
consider the following backward
stochastic heat equation:
\bel{bshe} \lt\{\!\!\!
\begin{array}{ll}
\ds d z(t,\!x) = \lt(-\D z(t,\!x)\! +F(t,x,z(t,\!x),Z(t,\!x) ) \rt) d t\!+\!Z(t,\!x) dW(t), &  (t,x) \in [0,T)\times D,\\
\ns\ds z(t,x) =0,    & (t,x) \in [0,T)\times \partial D,\\
\ns\ds z(T,x)=z_T(x), & x \in  D,
\end{array}
\rt.
\ee
where $D \subset \dbR^\ell$ (for some $\ell\in \dbN$) is a bounded domain
with a $C^2$ boundary $\pa D$.

\ss

Let us first introduce the following
two assumptions:

\ss

\textbf{(A1)} $\mathbf{F}$ is the natural filtration generated by $\{W(t)\}_{t\in[0,T]}$, $z_T\in L^2_{\mf_T}(\O;H_0^1(D))$, and $F:[0,T]\times D \times \dbR \times\dbR \to \dbR$ is $\frac{1}{2}$-H\"{o}lder continuous with respect to $t$ and Lipchitz continuous with respect to $z$ and $Z$. Moreover, $F(\cd,\cd,0,0)\in L^2(0,T;H_0^1(D))$.

\ss

\textbf{(A2)} For the (equal time-interval) partition $\pi:\,0=t_0<t_1<\cds<t_N=T$ with $N\in \dbN$, $\t=\frac T N$ and $t_i=i\t $ (for $i=0,1,\cds, N$),
there exists a constant $\cC>0$ such that the correction term $Z$ in the equation \rf{bshe} satisfies 
\bel{w909e1}
\sum_{n=0}^{N-1}\me\int_{t_n}^{t_{n+1}}\lt(\|Z(t)-Z(t_{n+1})\|_{L^2(D)}^2+\|Z(t)-Z(t_n)\|_{L^2(D)}^2\rt)dt\leq \cC \t.
\ee

\ss

\br{w909r4} 

{\rm(1).} The filtration $\mathbf{F}$ can be generalized to quasi-left continuous one, and $F(\cd,\cd,\cd,\cd)$ can be
stochastic. But we do not consider these general case to avoid  technical complexity.  

{\rm(2).} In this chapter, we adopt the equipartition for simplicity. Our method  can also be applied to quasi-uniform partition, i.e.
 $$\ds\max_{k=0,1,\cds,N-1}\{t_{k+1}-t_k\}\leq \cC\min_{k=0,1,\cds,N-1}\{t_{k+1}-t_k\}.$$

{\rm(3).} In \cite{ZhangJF04}, the estimate in the form of \rf{w909e1} is called the
$L^2$-regularity of $Z(\cd)$. Under
some suitable conditions, \rf{w909e1}
can be guaranteed \citep[see, e.g.,][]{Wang16,Wang20}. When
$Z(\cd)$ appears in the drift term,
\rf{w909e1} is crucial to prove rates of
convergence for temporal
discretization. If $Z(\cd)$ does not
appear in the drift term, when proving
the rate for $z(\cd)$, we do not need the condition
\rf{w909e1}  \citep[see][]{Prohl-Wang20}.
\er

We shall use the following notations: For $ k=0,1,\cds, N-1$, $t\in
[t_k,t_{k+1}),$
\bel{w9a2}
 \D_{k+1}W= W(t_{k+1})-W(t_k),\;\;
 \ee
\bel{w1e1}
\nu(t)= t_k,\;\;\mu(t)= t_{k+1},\;\; \pi(t)= k,
\ee
and
$$\nu(T)=\mu(T)= T,\;\;\pi(T) = N.$$
%

 
The numerical method (based on the
transposition solutions to  backward stochastic differential equations) given
in \cite{Wang-Wang-Zhang20} \citep[see also][]{Lv-Zhang13, Wang-Zhang11}
can be regarded as a Galerkin method, and hence we call it the 
finite transposition method. Let us propose below the
outline of the finite
transposition method to solve the equation \rf{bshe}:
\begin{enumerate}[(1)]

\item Determine  a finite dimensional subspace $\dbS $
 of $L^2(D)$.

\item  
Choose a {\em finite transposition space $\dbH=\mbox{\rm span}\{e_i\}$,} which is a finite dimensional subspace
 of $L_{\dbF}^2(0,T;\dbS)$ (in analogy with the finite element space of the finite element method in solving partial differential equations). 

\item Introduce the following variational equation:
\bel{w629e2}
\bal
&\me\langle \f(T),z_T\rangle_{L^2(D)}\\
&=\me\int_{0}^T  \[ \langle \f(t), F(\mu(t),\cd, z(t),Z(t))\ran_{L^2(D)}+ \langle
v_1(t),z(t)\rangle_{L^2(D)} +\langle v_2(t),
Z(t)\rangle_{L^2(D)} \] dt,
\eal
\ee
where 
$$ 
\f (t)=\int_0^{\mu(t)} \lt(\D\f(\nu(s))+ v_1(s)\rt)ds+\int_0^{\nu(t)} v_2(s)dW(s),
$$ 
with $v_1(\cd)$ and $v_2(\cd)$ being stochastic processes in
suitable finite dimensional  subspaces $\dbH_1\subset L^1_{\dbF}(t,T;L^2(\Omega;\dbS))$ and $\dbH_2\subset L_{\dbF}^2(0,T;\dbS)$, respectively, and $\nu(\cd)$ and $\mu(\cd)$ being piecewise constant  functions defined in \rf{w1e1}.
Based on this variational
equation, one can prove the existence and
uniqueness of  approximate solution  $(z(\cd),Z(\cd))$, called the {\em finite transposition solution}, to the
equation \eqref{bshe} in the form
\beq
z(\cd)=\sum_{i=1}^{\mbox{dim}(\dbH_1)} \a_i e_i,\q Z(\cd)=\sum_{i=1}^{\mbox{dim}(\dbH_2)} \b_i e_i.
\eeq

\item Find coefficients $\a_i,\, \b_i$ of the  {\em finite transposition solution}  via the variational equation \rf{w629e2},
and prove the rate of convergence.


\end{enumerate}

\begin{remark}
There is a useful method --- the stochastic finite
element method --- to solve partial differential equations with random parameter \citep[see, e.g.,][]{Ghanem-Spanos91}. A main ingredient  in this method is the
orthogonal expansions, such as
polynomial chaos expansion,
Karhunen-Lo\`eve expansion, etc.
Let $\{{\bf e}_i\}_{i=1}^\infty$ be an
orthogonal basis for
$L^2_{\cF_T}(\O;L^2(D))$. Then for any
stochastic process $X(\cdot)$, there
exists an expansion
\beq
X(\cdot)=\sum\limits_{i=1}^\infty x_{i}(\cd)
{\bf e}_i.
\eeq
Since the solutions $(z(\cd),Z(\cd))$
to backward stochastic evolution
equations are progressively measurable  stochastic
processes, people need to find the
orthogonal basis of
$L^2_\dbF(0,T;L^2(D))$. However, it
seems that there is no simple explicit
orthogonal basis for such a Hilbert
space.  Note that  $L^2_\dbF(0,T;L^2(D))$ is much
more complicated than
$L^2_{\cF_T}(\O;L^2(D))$, and hence
the classical
``stochastic finite element method"
does not work for our problem.

\end{remark}

Now we give details for the above
outline of the finite
transposition method.

\ss

Firstly, we determine a finite
dimensional subspace of $L^2(D)$ based
on the Galerkin method. Define
\bel{w923e1} A:\, D(A)=H_0^1(D)\cap
H^2(D)\to L^2(D),\q A\f=\D \f,
\q\forall \f\in D(A). \ee
 Let
$\{(-\l_i,\phi_i)\}_{i=1}^{\infty}$ be
the sequence of eigenvalues and
eigenfunctions of operator  $A$,  such
that $\|\phi_i\|_{L^2(D)} =1$ for
$i\in\dbN$. Then
$\{\phi_i\}_{i=1}^{\infty}$ constitutes
an orthonormal basis of $L^2(D)$ and
$\{\phi_i\}_{i=1}^{\infty}$ is also an
orthogonal basis of $H_0^1(D)$. 
Take
the subspace
\bel{w114e1}
\dbS_m=\mbox{span}\{\phi_1,\cdots,\phi_m\},\qq m\in \dbN.
\ee

Secondly, we construct a {\em
finite transposition space} by the
Wiener chaos. To this end, we
review some notations and results on
Wiener chaos \citep[see][for more details]{Nualart06}.
Recall that $\cF_{T}=\sigma\{W(t);0\leq
t\leq T\}$. Define the It\^o isometry
$\mathbb{W}:L^2(0,T)\to
L^2_{\cF_{T}}(\O)$ by
\beq \mathbb{W}(h)=\int_0^{T}h(t)
dW(t). \eeq
%
%
For $x\in\dbR$, let
\beq  H_n(x)= \left\{\!\!\!
\begin{array}{ll}
\ds {\frac{(-1)^n} {n!}}e^{\frac{x^2} {2}}\frac{d^n}{dx^n}(e^{-{\frac{x^2} {2}}}), & n>0,\\
\ns\ds 1, & n=0
\end{array}
\right. \eeq
be the $n$-th Hermite polynomial.
Denote by $\cI$ the set of all
sequences
$\alpha=(\alpha_1,\alpha_2,\cdots)$,
$\alpha_i\in\mathbb N_0$, such that all
the terms, except a finite number of
them, vanish. For $k\in\dbN$, let
$\cI(k)$ be the set of all sequences $\a=(\a_1,\a_2,\cds,\a_k),\,\a_i\in\dbN_0$. For all
$\alpha\in \cI \;\, (\resp\; \a\in\cI(k))$, write
\beq \alpha!\equiv
\prod\limits_{i=1}^\infty \alpha_i! \;\,
\Big(\resp\; \prod\limits_{i=1}^k
\alpha_i!\Big),\quad \mbox{and} \quad
|\alpha|\equiv \sum\limits_{i=1}^\infty
\alpha_i \;\, \Big(\resp\; \sum\limits_{i=1}^k \alpha_i\Big),
\eeq 
and for $x=(x_1,x_2,\cds)\in \dbR^\infty$,
\beq H_\alpha(x) =
\prod\limits_{i=1}^\infty
H_{\alpha_i}(x_i),\;\; \a\in\cI\; \,
  \Big(\resp\; \prod\limits_{i=1}^k H_{\alpha_i}(x_i),\;\; \a\in\cI(k)\Big),
\eeq
which is called the generalized Hermite polynomial.
Also, we put
\beq
g_i(\cd)=
\frac{\chi_{[t_{i-1},t_i)}(\cd)}{\sqrt{\t
}},\q i=1,2,\cds,N. \eeq
Furthermore, for any $k=1,2,\cds,N$,
define the {\em Wiener chaos of order
$n$} in $H$ as follows: 
\beq
\cH_n(k;H) \deq \span\Big\{
\prod_{i=1}^k H_{\a_i}(\dbW(g_i))\phi
\;\Big|\; \a\in
\cI(k),\,|\a|=n,\,\phi\in H \Big\}.
\eeq
Set $\cH_n(0;H)=H$.
Take
$\mathcal{G}_0=\mathcal{B}(\dbR)$, and
for $k=1,2,\cds,N$, define 
\beq
\mathcal{G}_{t_k} \deq
\si\{\D_1W,\D_2W,\cds,\D_kW\}, 
\eeq
where $\D_kW$ is defined in
\eqref{w9a2}. 
By virtue
of the Wiener chaos of order $n$, we
have the following orthogonal
decomposition result.

\begin{theorem}[{\citet[Theorem 1.1.1]{Nualart06}}]\label{w9l2}
For any $N,m\in \dbN$, and $k=0,1,\cds,N$, it holds that
\beq
L^2_{\cG_{t_k}}(\O;\dbS_m)=\bigoplus\limits_{n=0}^\infty \cH_n(k;\dbS_m).
\eeq
\et

Based on Theorem \ref{w9l2}, we set
\beq
\begin{array}{c}\ds
\cH^M(k;\dbS_m) \deq \bigoplus\limits_{n=0}^M \cH_n(k;\dbS_m),\\
\ns\ds \cH^M_N(k;\dbS_m) \deq \span \lt\{
\chi_{[t_{k},t_{k+1})}(\cd)\xi
\;\big|\; \xi\in \cH^M(k;\dbS_m)\rt\},
\end{array}
\eeq
and the desired {\em finite transposition space} as follows:
\bel{w9d2}
\begin{array}{c}\ds
\dbH_{N,M}(\dbS_m) \deq \bigoplus\limits_{k=0}^{N-1} \cH^M_N(k;\dbS_m)\subset L^2_\dbF(0,T;\dbS_m).
\end{array}
\ee
For simplicity, we denote
$\cH_n(k;\dbR)$, $\cH^M(k;\dbR)$,
$\cH_N^M(k;\dbR)$ and $\dbH_{N,M}(\dbR)$
by $\cH_n(k)$, $\cH^M(k)$, $\cH_N^M$
and $\dbH_{N,M}$, respectively. 
We write $\big(\cH^M(k)\big)^{\otimes n}$ and $\dbH_{N,M}^{\otimes n}$
for the $n$-copies of $\cH^M(k)$ and $\dbH_{N,M}$, respectively.

\br{w531r1}
In Theorem \ref{w9l2}, we only list a decomposition of $L^2_{\cG_{t_k}}(\O;\dbS_m)$, which is a subspace of
$L^2_{\mf_T}(\O;L^2(D))$. For the latter space, by letting
\beq
\cH_n(L^2(D)) \deq \span\Big\{&  \prod_{i=1}^\infty H_{\a_i}(\dbW(g_i))\phi_j \;\Big| \; \a\in \cI,\,|\a|=n,\,
\{g_i\}_{i=1}^\infty \mbox{ is an orthonormal }\\
&\mbox{ basis of } L^2(0,T),
\{\phi_j\}_{j=1}^{\infty} \mbox{ is an orthonormal basis of } L^2(D)
 \Big\},
\eeq
we have
\beq
L^2_{\mf_T}(\O;L^2(D))=\bigoplus\limits_{n=0}^\infty
\cH_n(L^2(D)).
\eeq
For the projection operator 
\bel{w0121e1}
\G_M:
L^2_{\mf_T}(\O;L^2(D))\rightarrow
\bigoplus\limits_{n=0}^M
\cH_n(L^2(D)),
\ee 
it is easy to check
that 
$\G_M\big(L^2_{\cG_{T}}(\O;\dbS_m))=\cH^M(N;\dbS_m)$. This fact is
crucial in the sequel. \er

\br{w9r1}
By the construction of the finite transposition space $\dbH_{N,M}$, it is easy to see that the orthonormal
basis of $\cH^M(k)$ is
$$\{h_{k,i}\}_{i=1}^{M_k}= \bigg\{\sqrt{\a!}\prod_{i=1}^k H_{\a_i}(\dbW(g_i))\;\Big|\;\a\in \cI(k), |\a|\leq M \bigg\},$$
which is a finite dimensional subspace of $L^2_{\cG_{t_k}}(\O)$.
Note that $M_k$ depends on $M$. For example, it is easy to check that 
\beq
M_k=
\lt\{\!\!\!
\begin{array}{cl}
 \ds k+1,  & M=1,\\
 \ds \frac{(k+1)(k+2)}{2}, & M=2,\\
 \ds \frac{k^3+8k^2+19k+6}{6}, & M=3.
\end{array}
\rt.
\eeq
Hence the orthonormal basis of $\dbH_{N,M}$ is
$$
\{e_{k,i}(\cd)\}_{k,i}= \lt\{\chi_{[t_k,t_{k+1})}(\cd)\frac{1}{\sqrt{\t }}h_{k,i}\;\Big|\; 1\leq i\leq M_k, 0\leq k\leq N-1\rt\},
$$
which is also a finite dimensional
subspace of $L^2_\dbF(0,T)$.

By virtue of $\{h_{k,i}\}_{i=1}^{M_k}$,
$\{e_{k,i}(\cd)\}_{1\leq i\leq M_k,
0\leq k\leq N-1}$ and
$\{\phi_j\}_{j=1}^m$, we can obtain the
orthonormal basis of $\cH^M(k;\dbS_m)$
and $\dbH_{N,M}(\dbS_m)$. \er

\br{w909r2}
In this chapter, since we mainly focus on presenting the idea of handling  the difficulties caused by the stochastic setting,  we choose the finite
dimensional subspace of $L^2(D)$ given by \eqref{w114e1}. One can also choose that space as the one in other numerical methods for solving deterministic partial differential equations, such as the finite element space. 

For the  infinite dimensional space of $L_\dbF^2(0,T)$, subspace constructed by characteristic functions instead of
Wiener chaos is also a choice \citep[e.g.,][]{Dai-Zhang-Zou17}.
\er

\br{w902r5}
A Wiener chaos decomposition method is  proposed in \cite{Briand-Labart14}.
It should be pointed out that, although in both \cite{Briand-Labart14} and the present chapter,  some tools of
Wiener chaos expansion are used, the numerical method in \cite{Briand-Labart14} (based on Malliavin analysis) is
totally different from the finite transposition method in this chapter (based on the variational equation \eqref{w629e2}).
\er

\section{Finite transposition method for backward stochastic evolution equations}\label{FTS-2}

Based on the finite transposition space $\dbH_{N,M}(\dbS_n)$, we can provide the definition of the finite transposition
solution.

\begin{definition}\label{defi-of-ftm1}
For any $N,M,n\in\dbN$, a couple $(z(\cdot), Z(\cdot))\in
\dbH_{N,M}(\dbS_n) \times \dbH_{N,M-1}(\dbS_n)$ is called a finite transposition
solution to the equation \eqref{bshe}, if for any
$v_1(\cd)\in \dbH_{N,M}(\dbS_n)$ and $v_2(\cd)\in \dbH_{N,M-1}(\dbS_n)$, the following variational equation
holds
\bel{w903e1}
\bal
&\me\llan \f(T),z_{N,T}^\pi\rran_{L^2(D)}\\
&=\me\int_{0}^T  \[ \llan \f(t),  F(\mu(t),\cd,z(t),Z(t))\rran_{L^2(D)}+ \llan v_1(t),z(t)\rran_{L^2(D)} +\llan v_2(t), Z(t)\rran_{L^2(D)} \] dt.
\eal
\ee
Here $z^\pi_{N,T}(\in \cH^M(N;\dbS_n))$
is an approximation of $z_T$, and $\f(\cd)$ solves the following forward equation:
\bel{w903e2}
\lt\{
\bal
&\f(t_k)=\int_0^{t_{k+1}}\lt(A\f(\nu(t))+v_1(t)\rt)dt+\int_{0}^{t_k}v_2(t)dW(t),\q k=1,2,\cds,N-1,\\
&\f(0)=0,\\
&\f(t)=\f(\nu(t)),\q t\in (0,T],
\eal \rt.
\ee
where $A$ is defined in \rf{w923e1}.

The  finite transposition method   for equation \eqref{bshe}  is the algorithm
to obtain the finite transposition solution by the variational equation \eqref{w903e1}.

\end{definition}

The rest of this section is divided  into two subsections. In the first one, we discretize
the original  equation with respect to space variables. In the second one, we discretize the obtained equation with respect to the time
variable.

\subsection{Spacial discretization of \eqref{bshe}}

Denote by $\Pi_n$ the orthognonal projection from $ L^2(D)$ to $\dbS_n$ and define $A_n$ by
$A_n=A|_{\dbS_n}$. 
The semi-discrete problem in space corresponding to the equation \eqref{bshe} is to find
a  pair $(z_n(\cd),Z_n(\cd))\in  L^2_{\dbF}(\O;C([0,T];\dbS_n))\times
L^2_{\dbF}(0,T;\dbS_n)$ solving the following equation: 
\bel{bsde1}
\left\{
\bal
&dz_n(t)=\big(-A_nz_n(t)+\Pi_n F(t, z_n(t),Z_n(t))\big)\mathrm dt+Z_n(t)\mathrm dW(t), \q t\in [0,T], \\
 \displaystyle
& z_n(T)=\Pi_n z_T.
\eal
\right.
\ee
The following result states the rate
of convergence for the solution of the equation \eqref{bsde1} to the one of the equation \eqref{bshe}.
\begin{theorem}[{\citet[Theorem 5.4]{WangY13}, \citet[Theorem 3.1]{Wang16}}]\label{convergence speed}%
Suppose that ${\rm (A1)}$ holds. Let
$(z(\cd),Z(\cd))$ and
$(z_n(\cd),Z_n(\cd))$ be solutions to
the equations \eqref{bshe} and
\eqref{bsde1}, respectively. Then the
following estimate holds \beq \bal
&\|z(\cd)-z_n(\cd)\|^2_{L^2_{\dbF}(\O;C([0,T];
L^2(D)))\cap L^2_{\dbF}(0,T; H_0^1(D))}
   +\|Z(\cd)-Z_n(\cd)\|^2_{L^2_{\dbF}(0,T;  L^2(D))}\\
&\leq \frac{\cC}{\l_{n+1}}\Big(\|z_T(\cd)\|^2_{L^2_{\mf_T}(\O; H_0^1(D))}
                        +\|F(\cd,\cd,0,0)\|^2_{L^2(0,T; H_0^1(D))}\Big).
\eal
\eeq

\end{theorem}


As we have explain in Remark \ref{w909r2}, to figure out the main part, $\dbS_n$ in the finite transposition space $\dbH_{N,M}(\dbS_n)$ is constructed by the first $n$
eigenfunctions of the operator $A$. We can also choose other finite element spaces $\dbV_h$, and construct
the finite transposition space $\dbH_{N,M}(\dbV_h)$. The readers are referred to \cite{Dunst-Prohl16, Prohl-Wang20}
for error estimates for finite-element based spacial discretization.


\subsection{Temporal discretization of  \eqref{bsde1}}


In this part, firstly we rewrite the
definition of the finite  transposition
solution and the finite transposition
method for the backward stochastic
differential equation \eqref{bsde1}.
Then we present the existence and
uniqueness of the finite transposition
solution. Finally, we give the rate
of convergence.

Since $\dbS_n=\mbox{span}\{\phi_1,\cdots,\phi_n\}$, and the solution $(z_n(\cd), Z_n(\cd))$ to the
equation  \eqref{bsde1} is
in the space $L^2_{\dbF}(\O;C([0,T];\dbS_n))\times L^2_{\dbF}(0,T;\dbS_n)$, we may take $(z_n(\cd), Z_n(\cd))$ to be the following form
\begin{equation}\label{expansion 1}
\begin{aligned}
z_n(\cdot)=\sum_{j=1}^n
a_{n,j}(\cdot)\phi_j,\q
Z_n(\cdot)=\sum_{j=1}^n
b_{n,j}(\cdot)\phi_j,
\end{aligned}
\end{equation}
where $a_{n,j}(\cdot)\in L^2_{\dbF}(\O;C([0,T] ))$ and $b_{n,j}(\cdot)\in L^2_{\dbF}(0,T )$, for  $j=1,2,\cdots,n$.
Set
\begin{equation*}
\begin{aligned}
a_n(\cdot)=\begin{pmatrix} a_{n,1}(\cdot)\\ a_{n,2}(\cdot)\\\vdots
\\ a_{n,n}(\cdot)
\end{pmatrix},
\q b_n(\cdot)=\begin{pmatrix} b_{n,1}(\cdot)\\ b_{n,2}(\cdot)\\\vdots
\\ b_{n,n}(\cdot)
\end{pmatrix},
\q \L_n=\begin{pmatrix}
-\l_1& 0  &\cdots &0\\
0   & -\l_2&\cdots&0\\
\vdots&\vdots&\ddots&\vdots\\
0&  0     &\cdots &-\l_n\\
\end{pmatrix},\\
\end{aligned}
\end{equation*}
and $$\ds F_n(\cdot, a_n(\cdot),
b_n(\cdot))=\begin{pmatrix}
\lan F(\cdot,\cd,z_n(\cdot),Z_n(\cdot)),\phi_1 \ran_{ L^2(D)}\\\lan F(\cdot,\cd,z_n(\cdot),Z_n(\cdot)),\phi_2 \ran_{ L^2(D)}\\\vdots\\\lan F(\cdot,\cd,z_n(\cdot),Z_n(\cdot)),\phi_n \ran_{ L^2(D)}\\
\end{pmatrix}.$$
Then $(a_n(\cd), b_n(\cd))$ solves the following backward stochastic differential equation:
\begin{equation}\label{bsden}
\left\{
\begin{aligned}
&da_n(t)=\big(-\L_n a_n(t)+F_n(t,a_n(t),b_n(t))\big) dt+b_n(t) dW(t), \q t\in [0,T], \\
& a_n(T)=\big({\lan z_T,\phi_1 \ran}_{ L^2(D)},\,{\lan z_T,\phi_2 \ran}_{ L^2(D)},\,\cdots,\,{\lan z_T,\phi_n \ran}_{ L^2(D)}\big)^\top.
\end{aligned}
\right.
\end{equation}

\begin{definition}\label{defi-of-ftm2}
For any $N,M,n\in\dbN$, a couple
$(a_{n,N,M}(\cdot),
b_{n,N,M-1}(\cdot))\in
\dbH_{N,M}^{\otimes n} \times
\dbH_{N,M-1}^{\otimes n}$ is called a
finite transposition solution to the
backward stochastic differential
equation \eqref{bsden} if for any
$v_1(\cd)\in \dbH_{N,M}^{\otimes
n}$ and $v_2(\cd)\in \dbH_{N,M-1}^{\otimes
n}$, the following variational equation
holds 
\bel{w704e3} \bal \me\llan
x(T),a_{n,T}^\pi\rran
=&\,\me\int_{0}^T  \[ \llan x(t),  F_n(\mu(t),a_{n,N,M}(t),b_{n,N,M-1}(t))\rran\\
&\qq\q+ \llan v_1(t),a_{n,N,M}(t)\rran
+\llan v_2(t), b_{n,N,M-1}(t)\rran \] dt,
\eal \ee where $x(\cd)$ is given by
\bel{w704e4} \lt\{ \bal
&x(t_k)=\int_0^{t_{k+1}}\lt(\L_n x(\nu(t))+v_1(t)\rt)dt+\int_{0}^{t_k}v_2(t)dW(t),\q k=1,\cds,N-1,\\
&x(0)=0,\\
&x(t)=x(\nu(t)),\q t\in (0,T],
\eal
\rt.
\ee
and $a^\pi_{n,T}\in \lt(\cH^M(N)\rt)^{\otimes n}$ is an approximation of $a_n(T)$.
\end{definition}

By choosing 
$
a_{n,T}^\pi=\big({\lan z_{N,T}^\pi,\phi_1 \ran}_{ L^2(D)},\,{\lan z_{N,T}^\pi,\phi_2 \ran}_{ L^2(D)},\,\cdots,\,{\lan z_{N,T}^\pi,\phi_n \ran}_{ L^2(D)}\big)^\top,
$
and combining Definitions
\ref{defi-of-ftm1} and
\ref{defi-of-ftm2}, 
we can see that
the pair of stochastic processes \linebreak
$\big(\sum_{j=1}^n
\lt(a_{n,N,M}(\cd)\rt)_j\phi_j,\sum_{j=1}^n
\lt(b_{n,N,M-1}(\cd)\rt)_j\phi_j \big)$
is just the finite
transposition solution to \rf{bshe}, where $(a_{n,N,M}(\cd))_j$ and $(b_{n,N,M-1}(\cd))_j$ are the $j$-th component of $a_{n,N,M}(\cd)$ and \linebreak$b_{n,N,M-1}(\cd)$ respectively.
Based on this, from now on,
we study \rf{w704e3} instead of
\rf{w903e1}.

By Remark \ref{w9r1}, we can denote the orthonormal basis of $\dbH_{N,M}^{\otimes n}$ by
$$\lt\{e_{k,i}^\ell (\cd) \; \big|\; 1\leq \ell\leq n,\, 1\leq i\leq M_k, 0\leq k\leq N-1\rt\}.$$
Here
\beq
e_{k,i}^\ell(\cd)=(0,\cds,0,e_{k,i}(\cd),0,\cds,0)^\top,
\eeq
the $\ell$-th component of which is $e_{k,i}(\cd)$. Similarly, we can define $h_{k,i}^\ell$.

\ss

The following result is addressed to
the existence and uniqueness  of the
finite transposition solution to
backward stochastic differential
equations. The readers are referred to
\cite{Wang-Wang-Zhang20} for the proof.

\begin{theorem}\label{w9t1}
Under the assumptions {\rm(A1)} and
{\rm(A2)}, for any $N,M,n\in\dbN$, the
backward stochastic differential
equation \eqref{bsden} admits a unique
finite transposition solution
$(a_{n,N,M}(\cd), b_{n,N,M-1}(\cd))\in
\dbH_{N,M}^{\otimes
n}\times\dbH_{N,M-1}^{\otimes n}$.
\end{theorem}

\br{w704r1} By Definition
\ref{defi-of-ftm2}, the variational
equation \eqref{w704e3} leads to an
implicit algorithm. We can also provide
other schemes based on the finite transposition method,
such as the following explicit algorithm:
\beq \bal \me\llan
x(T),a_{n,T}^\pi\rran
=&\me\int_{0}^T  \[ \llan x(t),  F_n(\mu(t),a_{n,N,M}(\mu(t)),b_{n,N,M-1}(\mu(t)))\rran\\
&\qq\q+ \llan v_1(t),a_{n,N,M}(t)\rran
+\llan v_2(t), b_{n,N,M-1}(t)\rran \] dt.
\eal \eeq
\er

By Remark \ref{w9r1}, we know that
$\{e_{k,i}(\cd)\,|\,1\leq i\leq M_k,\, 0\leq
k\leq N-1\}$ is an orthonormal basis of the finite
transposition space $\dbH_{N,M}$.
Hence, by Theorem \ref{w9t1}, we can
write
\bel{w23e1} 
\begin{array}{ll} \ds
a_{n,N,M}(\cd)=\sum_{\ell=1}^n\sum_{k=0}^{N-1}\sum_{i=1}^{M_k}\a_{k,i}^\ell
e_{k,i}^\ell (\cd),\\
\ns\ds\q \a_{k,i}^\ell \in\dbR \mbox{
for }i=1,\cds,M_k,\; k= 0,\cds,N-1,\;
\ell=1,\cds,n,
\\
\ns\ds
b_{n,N,M-1}(\cd)=\sum_{\ell=1}^n\sum_{k=0}^{N-1}\sum_{i=1}^{M_k'}\b_{k,i}^\ell
e_{k,i}^\ell (\cd),\\
\ns\ds\q \b_{k,i}^\ell \in\dbR \mbox{
for }i=1,\cds,M_k',\; k= 0,\cds,N-1,\;
\ell=1,\cds,n.
\end{array} 
\ee
In what follows, we would apply the
variational equation \eqref{w704e3} to
determine the coefficients of the
finite transposition solution
$(a_{n,N,M}(\cd),b_{n,N,M-1}(\cd))$.

By choosing $v_1(\cd)=e_{k,i}^\ell (\cd)$
and $v_2(\cd)=0$, we obtain
\beq  x(t)= \left\{\!\!
\begin{array}{ll}\ds
0, \q & t\in[0,t_k),\\
\ns\ds \L_0^{\pi(t)-k+1}\sqrt{\t } h_{k,i}^\ell , \q & t\in[t_k,T].\\
\end{array}
\right.
\eeq
Here 
\beq
\L_0=\lt(I_n-\L_n \t  \rt)^{-1}. 
\eeq 
By
virtue of \eqref{w704e3}, we have
\begin{eqnarray*}
&&\me\lan \L_0^{N-k}\sqrt{\t } h_{k,i}^\ell , a_{n,T}^\pi\ran\\
&&=\me\int_{t_k}^T\left\lan \L_0^{\pi(t)-k+1}\sqrt{\t } h_{k,i}^\ell , F_n\big(\mu(t),a_{n,N,M}(t),b_{n,N,M-1}(t)\big)\right\ran dt\\
&&\q+\frac{1}{\sqrt{\t }}\me\int_{t_k}^{t_{k+1}}\left\lan h_{k,i}^\ell , a_{n,N,M}(t)\right\ran dt\\
&&=\me\int_{t_k}^T\left\lan \L_0^{\pi(t)-k+1}\sqrt{\t } h_{k,i}^\ell ,  F_n\big(\mu(t),a_{n,N,M}(t),b_{n,N,M-1}(t)\big)\right\ran dt\\ 
&&\q+\frac{1}{\t }\me\int_{t_k}^{t_{k+1}}\bigg\lan h_{k,i}^\ell , \sum_{m=1}^n\sum_{l=0}^{N-1}\sum_{j=1}^{M_l}\a_{l,j}^m\chi_{[t_l,t_{l+1})}(t)h_{l,j}^m\bigg\ran dt\\ 
&&=\me\int_{t_k}^T\left\lan \L_0^{\pi(t)-k+1} \sqrt{\t } h_{k,i}^\ell ,  F_n\big(\mu(t),a_{n,N,M}(t),b_{n,N,M-1}(t)\big)\right\ran dt
+\a_{k,i}^\ell,
\end{eqnarray*}
which implies
\bel{w23e6}
\bal
\a_{k,i}^\ell =&\me\left\lan \L_0^{N-k}\sqrt{\t } h_{k,i}^\ell , a_{n,T}^\pi\right\ran\\
  &- \me\int_{t_k}^T\left\lan \L_0^{\pi(t)-k+1} \sqrt{\t } h_{k,i}^\ell ,  F_n\big(\mu(t),a_{n,N,M}(t),b_{n,N,M-1}(t)\big)\right\ran dt.
\eal
\ee

In the same vein, by taking $v_1=0,\,v_2=e_{k,i}^\ell $, we see that
\beq
x(t)=
\left\{\!\!
\begin{array}{ll}
0,\, & t\in[0,t_{k+1}),\\
\ns\ds \L_0^{\pi(t)-k}\frac{ \D_{k+1}W}{\sqrt{\t }} h_{k,i}^\ell ,\q &t\in[t_{k+1},T],\\
\end{array}
\right. \eeq
and
\begin{eqnarray*}
&&\me\Big\lan \L_0^{N-k-1} \frac{\D_{k+1}W}{\sqrt{\t }} h_{k,i}^\ell , a_{n,T}^\pi\Big\ran\\
&&=\me\int_{t_{k+1}}^T\Big\lan \L_0^{\pi(t)-k}\frac{\D_{k+1}W}{\sqrt{\t }} h_{k,i}^\ell , F_n\big(\mu(t),a_{n,N,M}(t),b_{n,N,M-1}(t)\big)\Big\ran dt\\
&& \q  +\frac{1}{\sqrt{\t }}\me\int_{t_k}^{t_{k+1}}\Big\lan h_{k,i}^\ell , b_{n,N,M-1}(t)\Big\ran dt\\
&&=\me\int_{t_{k+1}}^T\Big\lan \L_0^{\pi(t)-k} \frac{\D_{k+1}W}{\sqrt{\t }} h_{k,i}^\ell ,  F_n\big(\mu(t),a_{n,N,M}(t),b_{n,N,M-1}(t)\big)\Big\ran dt\\
&& \q +\frac{1}{\t }\me\int_{t_k}^{t_{k+1}}\Big\lan h_{k,i}^\ell , \sum_{m=1}^n \sum_{l=0}^{N-1}\sum_{j=1}^{M_l'}\b_{l,j}^m\chi_{[t_l,t_{l+1})}(t)h_{l,j}^m\Big\ran dt\\
&&=\me\int_{t_{k+1}}^T\Big\lan
\L_0^{\pi(t)-k}
\frac{\D_{k+1}W}{\sqrt{\t }}
h_{k,i}^\ell ,
F_n\big(\mu(t),a_{n,N,M}(t),b_{n,N,M-1}(t)\big)\Big\ran
dt +\b_{k,i}^\ell.
\end{eqnarray*}
It follows that
\bel{w23e9}
\bal
\b_{k,i}^\ell =&\, \me\Big\lan \L_0^{N-k-1}\frac{\D_{k+1}W}{\sqrt{\t }} h_{k,i}^\ell , a_{n,T}^\pi\Big\ran\\
\ns\ds&-\me\int_{t_{k+1}}^T\Big\lan
\L_0^{\pi(t)-k}\frac{\D_{k+1}W}{\sqrt{\t
}} h_{k,i}^\ell,
F_n\big(\mu(t),a_{n,N,M}(t),b_{n,N,M-1}(t)\big)\Big\ran
dt. \eal \ee

By virtue of  \eqref{w23e6} and
\eqref{w23e9}, we can obtain the finite transposition
solution $(a_{n,N,M}(\cd) ,$ $
b_{n,N,M-1}(\cd))$ to the backward
stochastic differential equation \rf{bsden}. Applying Wiener
chaos expansion, Theorem \ref{w23t1} given below
proposes the convergence rate of the
finite transposition method, and also
shows the relationship between the
Euler method
 and the finite transposition one.
The following space is needed in the
proof: 
\beq
D\big((-A)^{3/2}\big)=\Big\{\f=\sum_{i=1}^\infty
 \varrho_i \phi_i\; \Big|\; \sum_{i=1}^\infty \,
|\varrho_i|^2\l_i^3<\infty\Big\}.
\eeq

\bt{w23t1} 
Assume {\rm
(A1)} and {\rm (A2)}, and suppose that
for given $N,M,n\in\dbN$, $(a_n(\cd),b_n(\cd))$ $\in L^2_\dbF(\O;C([0,T];\dbR^n))\times L^2_\dbF(0,T;\dbR^n)$ and $(a_{n,N,M}(\cd),b_{n,N,M-1}(\cd))\in
\dbH_{N,M}^{\otimes
n}\times\dbH_{N,M-1}^{\otimes n}$
are adapted solution and finite transposition solution  to
the backward stochastic differential
equation \eqref{bsden}, respectively.  Then,
\begin{eqnarray}\nonumber
&&a_{n,N,M}(\cd)=\sum_{k=0}^{N-1}\chi_{[t_k,t_{k+1})}(\cd)\bigg[\me\lt(\L_0^{N-k}a_{n,T}^\pi\;\Big|\;\mf_{t_k}\rt)\\ \label{w23e10a}
&&\qq\qq-\me\Big(\int_{t_k}^T
\L_0^{\pi(t)-k+1}\G_M F_n(\mu(t),a_{n,N,M}(t),b_{n,N,M-1}(t)) dt \;\Big|\;\mf_{t_k}\Big)\bigg], \\
\label{w23e10b}
&&b_{n,N,M-1}(\cd)=\sum_{k=0}^{N-1}\chi_{[t_k,t_{k+1})}(\cd)\me\Big( \frac{\D_{k+1}W}{\t } a_{n,N,M}(t_{k+1})\;\Big|\;\mf_{t_k}\Big),
\end{eqnarray}
where $\L_0=\lt(I_n+\L_n \t \rt)^{-1}$. 
Furthermore, the following rate of convergence holds true:
\bel{w23e11}
\bal
&\sup_{0\leq t\leq T} \me|a_{n,N,M}(t)-a_n(t)|^2+\me\int_0^T|b_{n,N,M-1}(t)-b_n(t)|^2dt\\
&\leq \cC\Big[\me|a_n(T)-\bar a_n(T)|^2+\l_n^2 \t+\me\int_{0}^{T}\lt|(I_n-\G_M)F_n(s,a_n(s),b_n(s))\rt|^2 ds \Big],
\eal
\ee
where $\G_M$ is defined in \rf{w0121e1}.
Moreover, if $F(\cd,\cd,\cd,\cd)$ is linear with respect to the last two components, and $z_T\in L^2_{\mf_T}\lt(\O;D\big((-A)^{3/2}\big)\rt), F(\cd,\cd,0,0)\in L^2(0,T;H_0^1(D)\cap H^2(D))$, then it holds that
\bel{w0121e2}
\bal
\sup_{0\leq t\leq T} \me|a_{n,N,M}(t)-a_n(t)|^2+\me\int_0^T|b_{n,N,M-1}(t)-b_n(t)|^2dt
\leq \cC\Big[\me|a_n(T)-\bar a_n(T)|^2+ \t \Big].
\eal
\ee
 
\et

\begin{proof} 
The proof is long, and we carry out it by the following three steps.

\ms

{\bf Step 1.} In this step, we prove
\eqref{w23e10a}. By
\eqref{w23e1} and \eqref{w23e6}, we
know that
\beq
\bal
a_{n,N,M}(\cd)
=&\sum_{\ell=1}^n\sum_{k=0}^{N-1}\sum_{i=1}^{M_k}\bigg[ \me \left\lan \L_0^{N-k}\sqrt{\t } h_{k,i}^\ell , a_{n,T}^\pi\right\ran\\
&\q -\me\int_{t_k}^T\lt\lan
\L_0^{\pi(t)-k+1}\sqrt{\t }
h_{k,i}^\ell ,
F_n(\mu(t),a_{n,N,M}(t),b_{n,N,M-1}(t))\rt\ran
dt\bigg]\times e_{k,i}^\ell (\cd). \eal
\eeq
For the first term on the right side of the above equality, we have that
\bel{w23e13}
\bal
\me\lt\lan \L_0^{N-k} \sqrt{\t } h_{k,i}^\ell , a_{n,T}^\pi\rt\ran e_{k,i}^\ell (\cd)
&=\me\lt\lan  \L_0^{N-k}\sqrt{\t } h_{k,i}^\ell , a_{n,T}^\pi\rt\ran \chi_{[t_k,t_{k+1})}(\cd)\frac{1}{\sqrt{\t }}h_{k,i}^\ell \\
&=\chi_{[t_k,t_{k+1})}(\cd)\me\lt\lan
h_{k,i}^\ell ,
\L_0^{N-k}a_{n,T}^\pi\rt\ran
h_{k,i}^\ell . \eal \ee By the
definition of
$L^2_{\cG_{t_k}}(\O;\dbR^n),\,
\cH^M(k;\dbR^n)$ and Remark \ref{w9r1},
we can extend the orthonormal basis
$\{h_{k,i}^\ell \,|\, 1\leq \ell\leq n,\,
1\leq i\leq M_k\}$ of
$\lt(\cH^M(k)\rt)^{\otimes n}$  to the
orthonormal basis $\{h_{k,i}\,|\,1\leq
k\leq n,\, 1\leq i\leq \infty\}$ of
$L^2_{\cF_{t_k}}(\O;\dbR^n)$. Since for any $\eta \in L^2_{\cG_T}(\O;\dbR^n)$, $\G_M \eta\in
\lt(\cH^M(N)\rt)^{\otimes n}$,
\bel{w23e14}  \me(\G_M
\eta\;|\;\mf_{t_k})\in
\lt(\cH^M(N)\rt)^{\otimes n}\cap
L^2_{\cG_{t_k}}(\O;\dbR^n)=
\lt(\cH^M(k)\rt)^{\otimes n}. \ee 
By
\rf{w23e13}, \eqref{w23e14} and the
fact that $\{h_{k,i}^\ell \;|\; 1\leq
\ell\leq n,\,1\leq i\leq M_k\}$ is an
orthonormal basis of
$\lt(\cH^M(k)\rt)^{\otimes n}$ and
$a_{n,T}^\pi\in
\lt(\cH^M(N)\rt)^{\otimes n}$, we see
that
\bel{w23e16}
\bal
\sum_{\ell=1}^n \sum_{k=0}^{N-1}\sum_{i=1}^{M_k} \me\lt\lan  \L_0^{N-k} \sqrt{\t } h_{ki}^\ell , a_{n,T}^\pi\rt\ran e_{k,i}^\ell (\cd)=\sum_{k=0}^{N-1}\chi_{[t_k,t_{k+1})}(\cd)  \L_0^{N-k}\me\lt( a_{n,T}^\pi\;|\;\mf_{t_k}\rt).
\eal
\ee

Similarly,
\begin{eqnarray}\label{w23e17}
&&\sum_{\ell=1}^n\sum_{k=0}^{N-1}\sum_{i=1}^{M_k} \me\int_{t_k}^T\lt\lan  \L_0^{\pi(t)-k+1} \sqrt{\t } h_{k,i}^\ell ,  F_n\big(\mu(t),a_{n,N,M}(t),b_{n,N,M-1}(t)\big)\rt\ran dt\, e_{k,i}^\ell (\cd) \nonumber\\
&&=\sum_{\ell=1}^n\sum_{k=0}^{N-1}\chi_{[t_k,t_{k+1})}(\cd)\sum_{i=1}^{M_k}\me\int_{t_k}^T\lt\lan  h_{k,i}^\ell ,  \L_0^{\pi(t)-k+1} \G_M F_n\big(\mu(t),a_{n,N,M}(t),b_{n,N,M-1}(t)\big)\rt\ran dt\, h_{k,i}^\ell \nonumber \\
&&=\sum_{k=0}^{N-1}\chi_{[t_k,t_{k+1})}(\cd) \me\lt(\int_{t_k}^T  \L_0^{\pi(t)-k+1}\G_M F_n\big(\mu(t),a_{n,N,M}(t),b_{n,N,M-1}(t)\big) dt \;\Big|\;\mf_{t_k}\rt).
\end{eqnarray}
Therefore, by \eqref{w23e13}, \eqref{w23e16} and \eqref{w23e17}, we have \eqref{w23e10a}.

\ms

{\bf Step 2.} In this step, we prove
\eqref{w23e10b}.  
Noting that $a_{n,T}^\pi\in \lt(\cH^M(N)\rt)^{\otimes n}$, by Remark \ref{w9r1}, we can get that
\beq
a_{n,T}^\pi=\sum_{\ell=1}^n\sum_{m=0}^M\sum_{|\a|=m}d^{\a,\ell} \sqrt{\a!} \prod_{i=1}^{N} H_{\a_i}^\ell (\dbW(g_i)).
\eeq
Recalling that $H_n(x)$ is the  Hermite polynomial, we have
\beq
xH_n(x)=(n+1)H_{n+1}(x)-H_{n-1}(x).
\eeq
Consequently,
\beq
\bal
\D_{k+1}W a_{n,T}^\pi=&\sum_{\ell=1}^n \sum_{m=0}^M\sum_{|\a|=m}d^{\a,\ell} \sqrt{\a!} \prod_{i=1}^{N} H_{\a_i}^\ell (\dbW(g_i))\sqrt{\t }H_1(\dbW(g_{k+1}))\\
=&\sqrt{\t
}\sum_{\ell=1}^n\sum_{m=0}^M\sum_{|\a|=m}d^{\a,\ell}
\sqrt{\a!} \prod_{i=1}^{N}
\[(\a_{k+1}+1)H_{\bar\a_i}^\ell
(\dbW(g_i))-H_{\hat\a_i}^\ell
(\dbW(g_i))\], \eal
\eeq
where $\bar\a=\a+\g_{k+1}$ and
$\hat\a=\a-\g_{k-1}$. Here the $k$-th
component of $\g_k(\in\cI)$ is $1$, and the
others are $0$. By the definition of
$\bar \a$,  we have
$\bar\a_{k+1}\geq 1$. Therefore, by noting
that $|\h \a|\leq M-1$, we arrive at
\begin{eqnarray*}
&&\me(\D_{k+1}W a_{n,T}^\pi\; |\; \mf_{t_k})\\
&& = \sqrt{\t }\sum_{\ell=1}^n\sum_{m=0}^M\sum_{|\a|=m}d^{\a,\ell} \sqrt{\a!}\me\Big( \prod_{i=1}^{N} \[(\a_{k+1}+1)H_{\bar\a_i}^\ell (\dbW(g_i))-H_{\hat\a_i}^\ell (\dbW(g_i))\]\;\Big| \;\mf_{t_k}\Big)\\
&& = \sqrt{\t }\sum_{\ell=1}^n\sum_{m=0}^M\sum_{|\a|=m}d^{\a,\ell} \sqrt{\a!}(\a_{k+1}+1)\prod_{i=1}^{k} H_{\bar\a_i}^\ell (\dbW(g_i))\me\Big(\prod_{i=k+1}^{N} H_{\bar\a_i}^\ell (\dbW(g_i)) \; \Big|\;\mf_{t_k}\Big)\\
&&\q-\sqrt{\t }\sum_{\ell=1}^n\sum_{m=0}^M\sum_{|\a|=m}d^{\a,\ell} \sqrt{\a!}\me\Big( \prod_{i=1}^{N} H_{\hat\a_i}^\ell (\dbW(g_i))\;\Big|\;\mf_{t_k}\Big)\\
&&= -\sqrt{\t }\sum_{\ell=1}^n\sum_{m=0}^M\sum_{|\a|=m}d^{\a,\ell} \sqrt{\a!}\me\Big( \prod_{i=1}^{N} H_{\hat\a_i}^\ell (\dbW(g_i))\;\Big|\;\mf_{t_k}\Big)\\
&&\in \lt(\cH^{M-1}(k)\rt)^{\otimes n}.
\end{eqnarray*}
Furthermore, we can deduce that
\bel{w23el1} \bal
&\sum_{\ell=1}^n\sum_{i=1}^{M_k'}\me\Big\lan \frac{\D_{k+1}W}{\sqrt{\t }} h_{k,i}^\ell , a_{n,T}^\pi\Big\ran h_{k,i}^\ell \\
&=\sum_{\ell=1}^n\sum_{i=1}^{M_k'}\me\Big\lan  h_{k,i}^\ell , \me\Big( \frac{\D_{k+1}W}{\sqrt{\t }} a_{n,T}^\pi \; \Big| \; \mf_{t_k}\Big)\Big\ran h_{k,i}^\ell \\
&=\me\Big( \frac{\D_{k+1}W}{\sqrt{\t }} a_{n,T}^\pi \;\Big| \;\mf_{t_k} \Big).
\eal
\ee

With the same procedure, we can get that
\begin{eqnarray}\label{w23el50}
&&\sum_{\ell=1}^n\sum_{k=0}^{N-1}\sum_{i=1}^{M_k'} \me\int_{t_{k+1}}^T \Big\lan \frac{\D_{k+1}W}{\sqrt{\t }} h_{k, i}^\ell,   F_n\big(\mu(t),a_{n,N,M}(t),b_{n,N,M-1}(t)\big) \Big \ran dt\, e_{k,i}^\ell (\cd)\nonumber\\
&&= \sum_{k=0}^{N-1}\chi_{[t_k,t_{k+1})}(\cd)\sum_{\ell=1}^n\sum_{i=1}^{M_k'} \me \Big \lan h_{k,i}^\ell , \int_{t_{k+1}}^T \frac{\D_{k+1}W}{\t } \G_M F_n\big(\mu(t),a_{n,N,M}(t),b_{n,N,M-1}(t)\big) dt \Big\ran h_{k,i}^\ell \nonumber\\
&&= \sum_{k=0}^{N-1}\chi_{[t_k,t_{k+1})}(\cd) \me \Big(\int_{t_{k+1}}^T \frac{\D_{k+1}W}{\t }  \G_M F_n\big(\mu(t),a_{n,N,M}(t),b_{n,N,M-1}(t)\big)dt \;\Big|\;\mf_{t_k}\Big).
\end{eqnarray}
%
Combining \eqref{w23e1}, \eqref{w23e9},
\eqref{w23el1} and \eqref{w23el50}, we
conclude that
\begin{eqnarray*}
&\3n b_{n,N,M-1}(\cd)\3n&= \sum_{\ell=1}^n \sum_{k=0}^{N-1}\sum_{i=1}^{M_k'}\Big[ \me\Big\lan  \L_0^{N-k-1}\frac{\D_{k+1}W}{\sqrt{\t }} h_{k,i}^\ell, a_{n,T}^\pi\Big\ran\\
&&\q-\me\int_{t_{k+1}}^T\Big\lan \L_0^{\pi(t)-k} \frac{\D_{k+1}W}{\sqrt{\t }} h_{k,i}^\ell, \G_M F_n(\mu(t),a_{n,N,M}(t),b_{n,N,M-1}(t))\Big\ran dt \Big]e_{k,i}^\ell (\cd)\\
&&= \sum_{k=0}^{N-1}\chi_{[t_k,t_{k+1})}(\cd)\me\Big( \frac{\D_{k+1}W}{\t }\Big[ \L_0^{N-k-1}a_{n,T}^\pi\\
&&\qq\qq - \int_{t_{k+1}}^T \L_0^{\pi(t)-k} \G_M F_n(\mu(t),a_{n,N,M}(t),b_{n,N,M-1}(t)) dt\Big] \;\Big|\;\mf_{t_k}\Big)\\
&&= \sum_{k=0}^{N-1}\chi_{[t_k,t_{k+1})}(\cd)\me\Big( \frac{\D_{k+1}W}{\t } a_{n,N,M}(t_{k+1})\;\Big|\;\mf_{t_k}\Big),
\end{eqnarray*}
which is \eqref{w23e10b}.

\ms
 
{\bf Step 3.} In this step, by means of \eqref{w23e10a} and \eqref{w23e10b}, we prove the error estimates for the finite transposition method, which is a 
slight variation of \citep[Theorem
4.2]{Wang16}. We provide a sketch for completeness.

Set $\bar a_n(t_k)= a_{n,N,M}(t_k)\in L^2_{\mf_{t_k}}(\O;\dbR^n)$, for $k=0,1,\cds,N$. By martingale representation theorem, there exists a
square integrable process $\bar b_n(\cd)$, such that 
$$\bar a_n(t_{k+1})=\me(\bar a_n(t_{k+1})\, |\, \mf_{t_k})+\int_{t_k}^{t_{k+1}}\bar b_n(s)dW(s), \q k=0,1,\cds,N-1.$$
%
For $t\in [t_k,t_{k+1}),\,k=0,1,\cds,N-1$, define 
\beq
\bar a_n (t)= \lt(I_n-\L_n \t  \rt)\bar a_n(t_k)+(t-t_k)\G_M F_n(t_{k+1},\bar a_n(t_k),b_{n,N,M-1}(t_k))+\int_{t_k}^t \bar b_n (s)dW(s).
\eeq
Subsequently, by noting that $b_{n,N,M-1}(t_k)=\frac 1\t\me\big(\int_{t_k}^{t_{k+1}}\bar b_n(s)ds\,\big| \, \mf_{t_k}\big)$, 
difference between $(a_n(\cd),b_n(\cd))$ and $(\bar a_n(\cd),\bar b_n(\cd))$ can be estimated as follows:
\begin{eqnarray}\label{w0119e2}
&&\me|\lt(I_n-\L_n \t  \rt)(a_n(t_k)-\bar a_n(t_k))|^2+\me\int_{t_k}^{t_{k+1}}|b_n(s)-\bar b_n(s)|^2ds\nonumber\\
&&=\me\Big| (a_n(t_k)-\bar a_n(t_k))-\int_{t_k}^{t_{k+1}}\lt(F_n(s,a_n(s),b_n(s))-\G_M F_n(t_{k+1}, \bar a_n(t_k), b_{n,N,M-1}(t_k)) \rt)ds\Big|^2\nonumber\\
&&\leq (1+\cC\t)\me|a_n(t_{k+1})-\bar a_n(t_{k+1})|^2 \nonumber \\
&&\q+\cC\Big\{\me\int_{t_k}^{t_{k+1}}\lt[\lt| \L_n(a_n(t_j)-a_n(s))\rt|^2+\lt|(I_n-\G_M)F_n(s,a_n(s),b_n(s))\rt|^2\rt] ds \nonumber\\
&&\qq\q+\me\int_{t_k}^{t_{k+1}} \lt| \G_M\lt(F_n(t_{k+1},\bar a_n(t_k), b_{n,N,M-1}(t_k))-F_n(s,a_n(s),b_n(s)) \rt)\rt|^2 ds \Big\}.
\end{eqnarray}
By Gronwall's inequality and assumption {\rm(A2)}, we get that
\bel{w0119e3}
\bal
&\max_{k=0,1,\cds,N}\me|a_n(t_k)-a_{n,N,M}(t_k)|^2=\max_{k=0,1,\cds,N}\me|a_n(t_k)-\bar a_n(t_k)|^2\\
&\leq \cC\Big[\me|a_n(T)-\bar a_n(T)|^2+\t+\t \me\int_{0}^{T}\l_n^2 \lt(\lt| \L_n a_n(s)\rt|^2+\lt| b_n(s)\rt|^2+\lt|F_n(t,0,0)\rt|^2 \rt)ds\\
&\qq+\me\int_{0}^{T}\lt|(I_n-\G_M)F_n(s,a_n(s),b_n(s))\rt|^2 ds \Big]\\
&\leq \cC\Big[\me|a_n(T)-\bar a_n(T)|^2+\l_n^2 \t+\me\int_{0}^{T}\lt|(I_n-\G_M)F_n(s,a_n(s),b_n(s))\rt|^2 ds \Big].
\eal
\ee
Summing  \rf{w0119e2} from $k=0$ to $N-1$, applying \rf{w0119e3} and the fact $b_{n,N,M-1}(t_k)=\frac 1\t\me\big(\int_{t_k}^{t_{k+1}}\bar b_n(s)ds\,\big|\,\mf_{t_k}\big)$, we can derive \rf{w23e11}.
In the similar vein, we can prove \rf{w0121e2}. That completes the proof.
 
\end{proof}

\br{w912r2} By formulas
\rf{w23e10a} and \rf{w23e10b}, for
a linear backward stochastic evolution
equation,  the transposition method is
just the Euler method under
$\cH^M(N;\dbS_n)$-valued approximations
of the terminal value. Nevertheless, for
nonlinear equations, the finite
transposition method is different from
the Euler method. For a linear
equation, since the terminal value
is approximated by
$\cH^M(N;\dbS_n)$-valued random
variables, thanks to the variational
equation \rf{w903e1}, there is no need
to calculate conditional expectations.

\er

\section{Numerical method for optimal controls}\label{Ap}

In this section, we present an
application of the finite transposition
method. To avoid technical complexity,
we consider a simple SLQ problem.  To
be specific, we consider a cost
functional
\beq
{\mathcal J}( y_0;u(\cd))=\frac 1 2 \me
\Big[\int_0^T\big( \| y(t)
\|_{L^2(D)}^2+ \| u(t)
\|^2_{L^2(D)}\big)dt +  \| y(T)
\|^2_{L^2(D)}\Big]
\eeq
subject to the (controlled forward)
stochastic heat equation with
an additive noise
\bel{w212e1}
 \lt\{
\bal
&dy(t) =\big(A y(t) +u(t) \big)dt+\si dW(t), \q t\in (0,T],\\
& y(0)=y_0.
\eal
\rt.
\ee
Here $A$ is defined in \rf{w923e1}, and
$y_0,\,\si\in H_0^1(D)\cap H^2(D)$. Let
us now state a SLQ problem as follows:

\ms

\no{\bf Problem (SLQ)}. Search for
$\bar u(\cd)\in L^2_\dbF(0,T;L^2(D))$,
such that \beq \cJ(y_0;\bar
u(\cd))=\inf_{u(\cd)\in
L^2_\dbF(0,T;L^2(D))}\cJ(y_0;u(\cd)).
\eeq


By \cite{Lv-Zhang15}, the solvability
of Problem (SLQ) is equivalent to the
solvability of  the following
forward-backward stochastic evolution
equation:
\bel{w212e3}
\lt\{\!\!\!
\begin{array}{ll}
\ds d \bar y(t)
=\big(A \bar y (t)+ \bar u(t) \big)dt
+\si dW(t), \qq & t\in (0,T],\\
\ns \ds d z(t) =\big(-A z(t) + \bar y(t) \big) dt+Z(t) dW(t), \q & t\in [0,T),\\
\ns \ds \bar y(0)=y_0 ,\qquad z(T)=-\bar y(T) , &
\end{array}
\rt.
\ee
with the condition \bel{w212e4}
 \bar u(t)-z(t)=0, \q \ae (t,\om)\in [0,T]\times \O.
\ee

In the following, we mainly propose the
discretization of
\rf{w212e3}-\rf{w212e4}. To be
specific, for the forward equation
\rf{w212e3}$_1$, we adopt the
time-implicit Galerkin method  \citep[see, e.g.,][]{Grecksch-Kloeden96}, while
for the backward one \rf{w212e3}$_2$,
we apply the finite transposition
method proposed in Section \ref{FTS-2};
see also Algorithm \ref{alg2} below.

Since \rf{w212e3}-\rf{w212e4} is
coupled, how to obtain the convergence
of this discretization strategy? 
In
what follows, we  adopt optimal control
theory to deduce the convergence rates.
To do this, we discretize Problem (SLQ)
within two steps: firstly, we get the 
spacial semi-discretization (which is
referred to as Problem (SLQ)$_{S}$); secondly, we obtain a
spatio-temporal discretization (which
is referred to as Problem
(SLQ)$_{ST}$). Now, we present these two
discretizations.

\ms

\no{\bf Problem (SLQ)$_S$.}  For a
fixed $n\in\dbN$, minimize the
following cost functional over
$L^2_\dbF(0,T;$ $\dbS_n)$: \beq
\cJ(\Pi_n y_0;u(\cd))=\frac 1 2
\me\Big[\int_0^T\big(\|y_n(t)\|_{L^2(D)}^2+\|u_n(t)\|_{L^2(D)}^2\big)d
t+\|y_n(T)\|_{L^2(D)}^2\Big] \eeq
subject to the following
$\dbS_n$-valued stochastic differential
equation:
\bel{w223e2} \lt\{ \bal
&d y_n(t)=\big(A_n  y_n(t)+u_n(t)\big)d t+\Pi_n\si d W(t), \q t \in [0,T],\\
&y_n(0)=\Pi_n y_0=\sum_{k=1}^n\lan y_0,\phi_k\ran_{L^2(D)}\phi_k
. \eal \rt. \ee

By \cite{Yong-Zhou99}, the solvability
of Problem (SLQ)$_S$ is equivalent to
the solvability of the following
forward-backward stochastic
differential equation:
\bel{w908e1} 
\lt\{\!\!\! 
\begin{array}{ll} 
\ds d \bar
y_n(t)=\big(A_n  \bar y_n(t)+ \bar
u_n(t) \big)dt
+\Pi_n\si dW(t),\qquad &   t\in (0,T],\\
\ns\ds d z_n(t)=\big(-A_n  z_n(t)+ \bar y_n(t) \big) dt+Z_n(t) dW(t), \qquad  &   t\in [0,T),\\
\ns\ds \bar y_n(0)=\Pi_n y_0,\qquad
z_n(T)=-\bar y_n(T)  , &
\end{array} 
\rt. 
\ee
with the condition \bel{w908e2}
 \bar u_n(t)-z_n(t)=0, \q \ae (t,\om)\in [0,T]\times \O .
\ee

The following result is on the rate of
convergence  for spacial
discretization of Problem (SLQ).
\bt{rate-space}
Let $(\bar y(\cd),\,\bar
u(\cd))$ be the optimal pair of Problem
(SLQ) and $(\bar y_n(\cd),\,\bar
u_n(\cd))$ be the solution of Problem
(SLQ)$_S$ for $n\in \dbN$. Then, there
exists a constant $\cC$ such that
\begin{eqnarray*}
{\rm (i)} &&
\ds\me\big(\sup_{t\in[0,T]}\|\bar u(t)-\bar u_n(t)\|_{L^2(D)}^2\big)\leq \frac{\cC}{\l_{n+1}},\\
{\rm (ii)} && \ds\me\big(\sup_{t\in[0,T]}\|\bar
y(t)-\bar 
y_n(t)\|_{L^2(D)}^2\big)+\me\int_0^T\|\bar
y(t)-\bar y_n(t)\|_{H_0^1(D)}^2d t \leq
\frac{\cC}{\l_{n+1}}.\\
\end{eqnarray*}
%
\et

\begin{proof}

Applying It\^o's
formula to $\|\nabla \bar
y(\cd)\|_{L^2(D)}^2+\|\bar
y(\cd)\|_{L^2(D)}^2$, we obtain \beq
\sup_{t\in[0,T]}\me\|\bar
y(t)\|_{H_0^1(D)}^2\leq
\cC\Big[\|y_0\|_{H_0^1(D)}^2+\me\int_0^T\lt(\|\bar
u(t)\|_{L^2(D)}^2+\|\si\|_{L^2(D)}^2\rt)d
t\Big]. \eeq 
On the other hand, 
Theorem \ref{convergence speed} yields that
\beq
\bal
&\me\big(\sup_{t\in[0,T]}\|z(t)-z_n(t)\|_{L^2(D)}^2\big)
\leq \frac {\cC}{\l_{n+1}}\Big(\me\|\bar y(T)\|_{H_0^1(D)}^2+\me\int_0^T\|\bar y(t)\|_{H_0^1(D)}^2d t\Big).
\eal \eeq
%
The maximum conditions \eqref{w212e4}
and \eqref{w908e2}  and above two inequalities lead to assertion (i).

It\^o's formula to
$\|\bar
y(\cd)- \bar y_n(\cd)\|_{L^2(D)}^2$, 
Gronwall's inequality and assertion (i) lead to
\bel{w0123e1}
\bal
&\sup_{t\in[0,T]}\me\|\bar y(t)-\bar y_n(t)\|_{L^2(D)}^2+\me\int_0^T \|\nabla(\bar y(s)-\bar y_n (s)) \|^2_{L^2(D)} ds\\
&\leq \cC\Big[\|y_0-\Pi_n y_0\|_{L^2(D)}^2
+\me \int_0^T \big[ \|\bar u(s)-\bar u_n (s) \|^2_{L^2(D)}+ \|\si-\Pi_n\si\|^2_{L^2(D)}\big]ds\Big]\\
&\leq \frac {\cC}{\l_{n+1}}.
\eal
\ee
Applying It\^o's formula to $\|\bar
y(\cd)- \bar y_n(\cd)\|_{L^2(D)}^2$ again, then using Burkholder-Davis-Gundy inequality and \rf{w0123e1},
we can derive assertion (ii).
\end{proof}

Before presenting the spatio-temporal
discretization of Problem (SLQ), we
introduce the following two spaces:
\begin{eqnarray*}
\begin{array}{c}
\ds {\mathbb X}_{\t} \deq \lt\{x(\cd) \in L^2_\dbF\big(0,T;\dbS_n \big)\;\big|\; x(t)=x(t_k), \,\,   \forall t\in [t_k, t_{k+1}),  \,\, k=0,1,\cds, \, N-1\rt\} ,\\
\ns\ds {\mathbb U}_{\t}\deq \lt\{u(\cd) \in
L^2_\dbF\big(0,T;\dbS_n \big)\;\big|\;
u(t)=u(t_k),  \,\,    \forall t\in
[t_k, t_{k+1}),  \,\, k=0,1,\cds, \,
N-1\rt\} ,
\end{array}
\end{eqnarray*}
and for any $x(\cd)\in
\dbX_{\t}$ and $u(\cd)\in \dbU_{\t}$, 
\beq
\|x(\cd)\|_{\dbX_{\t}}\deq\Big(\t
\sum_{n=1}^N\me \|x(t_n)\|_{L^2(D)}
\Big)^{1/2},\qq
\|u(\cd)\|_{\dbU_{\t}}\deq\Big(\t
\sum_{n=0}^{N-1}\me \|u(t_n)\|_{L^2(D)}
\Big)^{1/2}. 
\eeq

\no{\bf Problem (SLQ)$_{ST}$.} For fixed $N,n\in\dbN$, minimize
the following  functional over
$\dbU_\t$: 
\beq
\bal
\cJ_{ST}(\Pi_n y_0;u_{nN}(\cd))=\frac 1
2\lt(\|y_{nN}(\cd)\|^2_{\dbX_{\t}}+\|u_{nN}(\cd)\|^2_{\dbU_{\t}}
\rt)+\frac 1 2\me\lt\| y_{nN}(T)
\rt\|_{L^2(D)}^2 , \eal \eeq 
where
$(y_{nN}(\cd),\, u_{nN}(\cd))$ satisfies  
\bel{w226e12} \lt\{ \bal
&y_{nN}(t_{k+1})-y_{nN}(t_k)= \tau
\big[A_n  y_{nN}(t_{k+1})+
u_{nN}(t_k) \big]\\
&\qq\qq\qq\qq\qq+\Pi_n\si \D_{k+1}W,\qq  k=0,1,\cds,N-1 ,\\
& y_{nN}(0)= \Pi_n y_0.
\eal \rt. \ee

The following result, which guarantees the
solvability of Problem (SLQ)$_{ST}$, is
a variant of \citep[Theorem
4.2]{Prohl-Wang20}.  
\bt{MP} Problem (SLQ)$_{ST}$ admits a
unique minimizer $(\bar y_{nN}(\cd), \bar
u_{nN}(\cd))\in {\mathbb X}_{\tau} \times
{\mathbb U}_{\tau}$, which satisfies the following coupled equation
for $0\leq k\leq N-1$:
\bel{w1212e3} \begin{cases}
\ds\bar y_{nN}(t_{k+1})=\L_0 \bar y_{nN}(t_k)+\t \L_0 \bar u_{nN}(t_k)+\L_0 \Pi_n\si \D_{k+1}W , \\
\ns\ds z_{nN}(t_k)=\L_0 \me\lt(z_{nN}(t_{k+1})-\t \bar y_{nN}(t_{k+1}) \big|\mf_{t_k}\rt), \\
\ns\ds \bar y_{nN}(0)=y_{nN}(0), \q
z_{nN}(T)=-\bar y_{nN}(T), \end{cases}
\ee
together with the maximum condition
\beq
\bar u_{nN}(t_k)-z_{nN}(t_k)=0 \qquad
k=0,1,\cds,N-1,\q\as\!\!,
\eeq
where $\L_0=\lt(I_n-\L_n \t  \rt)^{-1}.$

\et

Based on \cite{Malanowski82}, as well
as the regularity of the optimal pairs
$(\bar y_n(\cd),\, \bar u_n(\cd))$ and
$(\bar y_{nN}(\cd),$ $\bar
u_{nN}(\cd))$, we can prove the
following convergence rates \citep[e.g.,][Theorem
4.3]{Prohl-Wang20}.

\bt{rate3}
For any
$n,N\in\dbN$, assume that $(\bar
y_n(\cd),\, \bar u_n(\cd))$ and $(\bar
y_{nN}(\cd),\, \bar u_{nN}(\cd))$ are
the optimal pairs of Problem (SLQ)$_S$ and
(SLQ)$_{ST}$, respectively. Then, it
holds that
\begin{eqnarray*}
{\rm (i)} &&
\sum_{k=0}^{N-1}\me\int_{t_k}^{t_{k+1}}\lt\|\bar
u_n (t)-\bar
u_{nN}(t_k)\rt\|_{L^2(D)}^2 d t
     \leq  \cC \t ; \\
{\rm (ii)}&& \max_{0 \leq k \leq N}
\me\lt\|\bar y_n(t_k)-\bar
y_{nN}(t_k)\rt\|_{L^2(D)}^2
 +\t\me \sum_{k=1}^N   \lt\|\bar y_n(t_k)-\bar y_{nN}(t_k)\rt\|_{H_0^1(D)}^2
 \leq \cC \t.\\
\end{eqnarray*}
\et

\br{w920r1}

The equation \rf{w1212e3} is the
implicit Euler scheme for the equation
\rf{w908e1} and the Galerkin based implicit Euler method for \rf{w212e3}.
Notice that \rf{w212e3} and \rf{w908e1} are coupled forward-backward equations.
Numerical methods for these coupled equations and the convergence analysis is highly nontrivial. To deduce the rates of convergence, additional conditions such as sufficiently small $T$, or equations' weak coupling, are needed \citep[see, e.g.,][]{Bender-Zhang08}. In this section, we
show an idea to solve these equations. To be specific, if a coupled equation can be transferred
to an optimal control problem, one can try to utilize tools in optimal control theory.
\er

The equation \rf{w1212e3} is a coupled equation and contains conditional expectations. There are some methods to
deal with these conditional expectations. But for coupled forward-backward equation, even for that in finite dimensions,
numerical schemes are rare and the convergence analysis is far from complete \citep[see, e.g.,][]{Bender-Zhang08}.
Here we can adopt the finite transposition method and gradient method to solve \rf{w212e3}.

\begin{algorithm}
\caption{\small Solving Problem (SLQ) by Galerkin based finite transposition method and the gradient method} \label{alg2}

Fix $N,M,n\in\dbN$, and $\kappa > 0$. Let $\wt u^{(0)}\in \dbH_{N,M}(\dbS_n)$. For any $\ell \in {\mathbb N}$,
update $\wt u^{(\ell)} \in \dbH_{N,M}(\dbS_n)$ as follows:
\begin{enumerate}[1.]
\item Compute $\wt y^{(\ell)}\in \dbH_{N,M}(\dbS_n)$ via time-implicit  Galerkin method:
\beq
\left\{
\bal
&\wt y^{(\ell)}(t_{n+1})-\wt y^{(\ell)}(t_{n})=\t A_n \wt y^{(\ell)}(t_{n+1})+\Pi_n \si  \D_{n+1}W,\q n=0,1,\cds,N-1,\\
& \wt y^{(\ell)}(0)=\Pi_n y_0 .
\eal
\right.
\eeq
\item Utilize $\wt y^{(\ell)} \in \dbH_{N,M}(\dbS_n)$ in \rf{w212e3}$_2$ to compute $\wt z^{(\ell)}\in \dbH_{N,M}(\dbS_n)$ via the finite transposition method.

\item Compute the update $\wt u^{(\ell+1)} \in \dbH_{N,M}(\dbS_n)$ via
\beq
\wt u^{(\ell+1)}=\big(1-\frac 1 {\kappa}\big)\wt u^{(\ell)}+\frac 1 \kappa \wt z^{(\ell)}  \, .
\eeq
\end{enumerate}
\end{algorithm}

In Problem (SLQ), $\kappa$ can be taken as $1+T+T^2$. We can
prove that $\lt(\wt u^{(k)}(\cd),
\wt y^{(k)}(\cd)\rt)$ converges to
$\lt(\bar u(\cd), \bar
y(\cd)\rt)$ in $\dbX_\t\times
\dbU_\t$ as $k\to\infty$ \citep[e.g.,][Section
5]{Prohl-Wang20}.


{\small
\begin{thebibliography}{32}

\bibitem[{Al-Hussein(2009)}]{Al-Hussein09}
Al-Hussein, A., 2009. Backward stochastic partial differential equations driven
  by infinite-dimensional martingales and applications. Stochastics 81~(6),
  601--626.

\bibitem[{Bally and Pag{\`e}s(2003)}]{Bally-Pages03}
Bally, V., Pag{\`e}s, G., 2003. Error analysis of the optimal quantization
  algorithm for obstacle problems. Stochastic Process. Appl. 106~(1), 1--40.

\bibitem[{Bender and Denk(2007)}]{Bender-Denk07}
Bender, C., Denk, R., 2007. A forward scheme for backward {SDE}s. Stochastic
  Process. Appl. 117~(12), 1793--1812.

\bibitem[{Bender and Zhang(2008)}]{Bender-Zhang08}
Bender, C., Zhang, J., 2008. Time discretization and {M}arkovian iteration for
  coupled {FBSDE}s. Ann. Appl. Probab. 18~(1), 143--177.

\bibitem[{Bensoussan(1969)}]{Bensoussan69}
Bensoussan, A., 1969. Contr\^{o}le optimal stochastique de syst{\`e}me
  gouvern\'{e}s par des \'{e}quations aux d\'{e}riv\'{e}es partielles de type
  parabolique. Rend. Mat. (6) 2, 135--173.

\bibitem[{Bensoussan(1983)}]{Bensoussan83}
Bensoussan, A., 1983. Stochastic maximum principle for distributed parameter
  systems. J. Franklin Inst. 315~(5-6), 387--406.

\bibitem[{Bouchard and Touzi(2004)}]{Bouchard-Touzi04}
Bouchard, B., Touzi, N., 2004. Discrete-time approximation and {M}onte-{C}arlo
  simulation of backward stochastic differential equations. Stochastic Process.
  Appl. 111~(2), 175--206.


\bibitem[{Briand and Labart(2014)}]{Briand-Labart14}
Briand, P., Labart, C., 2014. Simulation of {BSDE}s by {W}iener chaos
  expansion. Ann. Appl. Probab. 24~(3), 1129--1171.


\bibitem[{Dai et~al.(2017)}]{Dai-Zhang-Zou17}
Dai, L., Zhang, Y., Zou, J., 2017. Numerical schemes for forward-backward
  stochastic differential equations using transposition solutions. preprint.

\bibitem[{Douglas et~al.(1996)}]{Douglas-Ma-Protter96}
Douglas, Jr., J., Ma, J., Protter, P., 1996. Numerical methods for
  forward-backward stochastic differential equations. Ann. Appl. Probab. 6~(3),
  940--968.

\bibitem[{Dunst and Prohl(2016)}]{Dunst-Prohl16}
Dunst, T., Prohl, A., 2016. The forward-backward stochastic heat equation:
  numerical analysis and simulation. SIAM J. Sci. Comput. 38~(5), A2725--A2755.

\bibitem[{E et~al.(2019)}]{E-Hutzenthaler-Jentzen-Kruse19}
E, W., Hutzenthaler, M., Jentzen, A., Kruse, T., 2019. On multilevel {P}icard
  numerical approximations for high-dimensional nonlinear parabolic partial
  differential equations and high-dimensional nonlinear backward stochastic
  differential equations. J. Sci. Comput. 79~(3), 1534--1571.

\bibitem[{Ghanem and Spanos(1991)}]{Ghanem-Spanos91}
Ghanem, R.G., Spanos, P.D., 1991. Stochastic Finite Elements: a
  Spectral Approach. Springer-Verlag, New York.

\bibitem[{Gobet et~al.(2005)}]{Gobet-Lemor-Warin05}
Gobet, E., Lemor, J.-P., Warin, X., 2005. A regression-based {M}onte {C}arlo
  method to solve backward stochastic differential equations. Ann. Appl.
  Probab. 15~(3), 2172--2202.

\bibitem[{Grecksch and Kloeden(1996)}]{Grecksch-Kloeden96}
Grecksch, W., Kloeden, P.~E., 1996. Time-discretised {G}alerkin approximations
  of parabolic stochastic {PDE}s. Bull. Austral. Math. Soc. 54~(1), 79--85.

\bibitem[{Hu et~al.(2011)}]{Hu-Nualart-Song11}
Hu, Y., Nualart, D., Song, X., 2011. Malliavin calculus for backward stochastic 
differential equations and application to numerical solutions. Ann. Appl. Probab. 21~(6),
  2379--2423.

\bibitem[{Hu and Peng(1991)}]{Hu-Peng91}
Hu, Y., Peng, S.~G., 1991. Adapted solution of a backward semilinear stochastic
  evolution equation. Stochastic Anal. Appl. 9~(4), 445--459.

\bibitem[{Kushner(1968)}]{Kushner68}
Kushner, H.~G., 1968. On the optimal control of a system governed by a linear parabolic 
equation with white noise inputs. SIAM J. Control 6, 596--614.

\bibitem[{Li and Zhou(2020)}]{Li-Zhou20}
Li, B., Zhou, Q., 2020. Discretization of a distributed optimal control problem with a stochastic parabolic equation
driven by multiplicative noise.
arXiv:2011.14587v2. 

\bibitem[{Lions and Magenes(1972)}]{Lions-Magenes72}
Lions, J.-L., Magenes, E., 1972. Non-homogeneous boundary value problems and
  applications. {V}ol. {I}. Springer-Verlag, New York-Heidelberg.

\bibitem[{L\"{u} et~al.(2012)}]{Lu-Yong-Zhang12} L\"{u}, Q., Yong, J., Zhang, X., 2012. Representation of Itô integrals by Lebesgue/Bochner integrals. J. Eur. Math. Soc. 14~(6), 1795--1823.

\bibitem[{L\"{u} et~al.(2018)}]{Lu-Yong-Zhang18} L\"{u}, Q., Yong, J., Zhang, X., 2018. Erratum to ``Representation of It\^o integrals by Lebesgue/Bochner integrals" (J. Eur. Math. Soc. 14, 1795–1823 (2012)) [MR2984588]. J. Eur. Math. Soc. 20~(1) (2018), 259--260.

\bibitem[{L\"{u} and Zhang(2013)}]{Lv-Zhang13}
L\"{u}, Q., Zhang, X., 2013. Well-posedness of backward stochastic differential
  equations with general filtration. J. Differential Equations 254~(8),
  3200--3227.
 

\bibitem[{L\"{u} and Zhang(2014)}]{Lv-Zhang14}
L\"{u}, Q., Zhang, X., 2014. General {P}ontryagin-type Stochastic Maximum
  Principle and Backward Stochastic Evolution Equations in Infinite Dimensions.
  SpringerBriefs in Mathematics. Springer, Cham.

\bibitem[{L\"{u} and Zhang(2015)}]{Lv-Zhang15}
L\"{u}, Q., Zhang, X., 2015. Transposition method for backward stochastic
  evolution equations revisited, and its application. Math. Control Relat.
  Fields 5~(3), 529--555.


\bibitem[{L\"{u} and Zhang(2020)}]{Lv-Zhang20}
L\"{u}, Q., Zhang, X., 2020. Mathematical Theory for Stochastic Distributed
  Parameter Control Systems. Springer (In press).


  
\bibitem[{Malanowski(1982)}]{Malanowski82}
Malanowski, K., 1982. Convergence of approximations vs. regularity of solutions
  for convex, control-constrained optimal-control problems. Appl. Math. Optim.
  8~(1), 69--95.

\bibitem[{Nualart(2006)}]{Nualart06}
Nualart, D., 2006. The {M}alliavin Calculus and Related Topics, second ed.
  Springer-Verlag, Berlin.

\bibitem[{Prohl and Wang(2020a)}]{Prohl-Wang20-2}
Prohl, A., Wang, Y., 2020a. Strong error estimates for a space-time discretization 
of the linear-quadratic control problem with the stochastic heat equation with linear noise.
arXiv:2012.04418, Submitted. 

\bibitem[{Prohl and Wang(2020b)}]{Prohl-Wang20}
Prohl, A., Wang, Y., 2020b. Strong rates of convergence for space-time
  discretization of the backward stochastic heat equation, and of a
  linear-quadratic control problem for the stochastic heat equation.
  https://na.uni-tuebingen.de/preprints.shtml, arXiv-2012.10117, Submitted. 

\bibitem[{Tzafestas and Nightingale(1968)}]{Tzafestas-Nightingale68}
Tzafestas, G., Nightingale, J.~M., 1968. Optimal control of a class of linear 
stochastic distributed- parameter systems. Proc. Inst. Electr. Engrs. 115, 1213--1220.

\bibitem[{Wang et~al.(2020)}]{Wang-Wang-Zhang20}
Wang, P., Wang, Y., Zhang, X., 2020. Error analysis of finite transposition
  method for solving backward stochastic differential equations. Preprint.

\bibitem[{Wang and Zhang(2011)}]{Wang-Zhang11}
Wang, P., Zhang, X., 2011. Numerical solutions of backward stochastic
  differential equations: a finite transposition method. C. R. Math. Acad. Sci.
  Paris 349~(15-16), 901--903.

\bibitem[{Wang(2013)}]{WangY13}
Wang, Y., 2013. Transposition Solutions of Backward Stochastic Differential
  Equations and Numerical Schemes. Ph.D. Thesis, Academy of Mathematics and
  Systems Science, Chinese Academy of Sciences.

\bibitem[{Wang(2016)}]{Wang16}
Wang, Y., 2016. A semidiscrete {G}alerkin scheme for backward stochastic
  parabolic differential equations. Math. Control Relat. Fields 6~(3),
  489--515.

\bibitem[{Wang(2020)}]{Wang20}
Wang, Y., 2020. {$L^2$}-regularity of solutions to linear backward stochastic
  heat equations, and a numerical application. J. Math. Anal. Appl. 486~(1),
  123870, 18.

\bibitem[{Yong and Zhou(1999)}]{Yong-Zhou99}
Yong, J., Zhou, X.~Y., 1999. Stochastic Controls: Hamiltonian Systems and HJB
  Equations. 
  Springer-Verlag, New York.

\bibitem[{Zhang(2004)}]{ZhangJF04}
Zhang, J., 2004. A numerical scheme for {BSDE}s. Ann. Appl. Probab. 14~(1),
  459--488.

\bibitem[{Zhao et~al.(2006)}]{Zhao-Chen-Peng06}
Zhao, W., Chen, L., Peng, S., 2006. A new kind of accurate numerical method for
  backward stochastic differential equations. SIAM J. Sci. Comput. 28~(4),
  1563--1581.

\bibitem[{Zhao et~al.(2014)}]{Zhao-Fu-Zhou14}
Zhao, W., Fu, Y., Zhou, T., 2014. New kinds of high-order multistep schemes for 
coupled forward backward stochastic differential equations. SIAM J. Sci. Comput. 36~(4),
  A1731--A1751.

\end{thebibliography}
}

\end{document}